\documentclass[reqno, oneside]{amsart}
\usepackage{style}

\begin{document}
\title[Duality theorems]{Duality theorems for coinvariant subspaces of $H^1$}
\author{R.~V.~Bessonov}

\address{St.Petersburg State University ({\normalfont \hbox{7-9}, Universitetskaya nab., 199034, St.Petersburg,  Russia}), St.Petersburg Department of Steklov Mathematical Institute of Russian Academy of Science ({\normalfont 27, Fon\-tan\-ka, 191023, St.Petersburg, Russia}), and School of Mathematical Sciences, Tel Aviv University ({\normalfont 69978, Tel Aviv, Israel})}
\email{bessonov@pdmi.ras.ru}
\thanks{This work is partially supported by RFBR grants 12-01-31492, 14-01-00748, by ISF grant~94/11, by JSC ``Gazprom Neft'', and by the Chebyshev Laboratory (Department of Mathematics and Mechanics, St. Petersburg State University) under RF Government grant 11.G34.31.0026}

\begin{abstract}
Let $\theta$ be an inner function satisfying the connected level set condition of B.~Cohn, and let $\Kl$ be the shift-coinvariant subspace of the Hardy space~$H^1$ generated by $\theta$. We describe the dual space to $\Kl$ in terms of a bounded mean oscillation with respect to the Clark measure $\sa$ of $\theta$. Namely, we prove that $(\Kl \cap zH^1)^* = \bmosa$. The result implies a two-sided estimate for the operator norm of a finite Hankel matrix of size $n\times n$ via \hbox{${\rm BMO}(\mu_{2n})$-norm} of its standard symbol, where $\mu_{2n}$ is the Haar measure on the group $\{\xi \in \C: \xi^{2n} = 1\}$.
\end{abstract}

\maketitle

\section{Introduction}\label{s1}
A bounded analytic function $\theta$ in the open unit disk $\D = \{z \in \C: |z| < 1\} 	$ is called inner if \hbox{$|\theta(z)| = 1$} for almost all points $z$ on the unit circle $\T$ in the sense of angular boundary values. With every inner function $\theta$ we associate the shift-coinvariant \cite{Nik86} subspace $\Kthp$ of the Hardy space $H^p$,
\begin{equation}\label{eq1}
\Kthp = H^p \cap \bar z \theta \ov{H^p}, \quad 1 \le p \le \infty.
\end{equation}
As usual, functions in $H^p$ are identified with their angular boundary values on the unit circle $\T$; formula \eqref{eq1} means that $f \in \Kthp$ if $f \in H^p$ and there is $g \in H^p$ such that $f(z) = \bar z \theta(z) \ov{g(z)}$ for almost all points $z \in \T$. An inner function $\theta$ is said to be {\it one-component} if its sublevel set $\Omega_\delta = \{z \in \D: |\theta(z)|< \delta\}$ is connected for a positive number $\delta < 1$. This class of inner functions was introduced by B.Cohn \cite{Co82} in 1982. It is very useful in studying   Carleson-type embeddings $\Kthp \into L^p(\mu)$ and Riesz bases of reproducing kernels in~$\Kthp$, see \cite{Co82,Co86,VT86,Al99,BD11,BBK11, AR1, AR2} for results and further references.   

\medskip

In this paper we describe the dual space to the space $\Kl$ generated by a 
one-component inner function $\theta$. Our main result is the following formula: 
\begin{equation}\label{eq36}
(\Kl \cap z H^1)^* = \bmosa,
\end{equation}
where $\sigma_\alpha$ denotes the Clark measure of the inner function $\theta$. Below we state this result formally and apply it to the boundedness problem for truncated Hankel operators. 

\subsection{Clark measures of one-component inner functions}
Let $\theta$ be a non-constant inner function in the open unit disk~$\D$. For each complex number $\alpha$ of unit modulus the function $\Re\bigl(\frac{\alpha + \theta}{\alpha - \theta}\bigr)$ is positive and harmonic in $\D$. Hence there exists the unique positive Borel measure $\sigma_\alpha$ supported on the unit circle $\T$  such that  
\begin{equation}\label{eq20}
\Re \frac{\alpha + \theta(z)}{\alpha - \theta(z)} = \int_{\T}\frac{1 - |z|^2}{|1 - \bar \xi z|^2}\,d\sigma_\alpha(\xi), \quad z \in \D.
\end{equation}
The measures $\{\sigma_\alpha\}_{|\alpha| = 1}$ are usually referred to as Clark measures of the inner function~$\theta$ due to seminal work \cite{Cl72} of D.~N.~Clark where their close connection to rank-one perturbations of singular unitary operators was discovered. For a modern exposition of this topic and subsequent results see survey  \cite{PoltSar06}.

Each Clark measure $\sigma_\alpha$ of an inner function $\theta$ is singular with respect to the Lebesgue measure on the unit circle $\T$. Conversely if $\mu$ is a finite positive Borel singular measure supported on $\T$ and $|\alpha| = 1$, then there exists the unique inner function~$\theta$ satisfying~\eqref{eq20} with $\sigma_\alpha = \mu$. Thus, there is one-to-one correspondence between inner functions in the unit disk $\D$ and singular measures on the unit circle~$\T$. It was unknown which singular measures on $\T$ correspond to the Clark measures of one-component inner functions. We fill this gap in Theorem \ref{t1} below. 

\medskip

For every Borel measure $\mu$ on the unit circle $\T$ denote by $a(\mu)$ the set of isolated atoms of $\mu$. Then the set $\rho(\mu) = \supp \mu \setminus a(\mu)$ consists of accumulating points in the support $\supp \mu$ of $\mu$. We will say that an atom $\xi \in a(\mu)$ has two neighbours in~$a(\mu)$ if there is an open arc $(\xi_-, \xi_+)$ of the unit circle $\T$ with endpoints~$\xi_\pm \in a(\mu)$ such that $\xi$ is the only point in $(\xi_-, \xi_+) \cap \supp\mu$. By $m$ we will denote the Lebesgue measure on~$\T$ normalized so that $m(\T) = 1$. 
\begin{Thm}\label{t1} Let $|\alpha| = 1$. 
The following conditions are necessary and sufficient for a Borel measure $\mu$  to be the Clark measure $\sigma_\alpha$ of a one-component inner function:
\begin{itemize}
\item[(a)] $\mu$ is a discrete measure on $\T$ with isolated atoms,  $m(\supp \mu) = 0$, every atom $\xi \in a(\mu)$ has two neighbours $\xi_\pm$ in $a(\mu)$, and every connected component of $\T\setminus \rho(\mu)$ contains  atoms of $\mu$; 
\item[(b)] $A_\mu |\xi-\xi_\pm| \le \mu\{\xi\} \le B_\mu |\xi-\xi_\pm|$ for all $\xi \in a(\mu)$ and some  $A_\mu> 0,$ $B_\mu < \infty$;
\item[(c)] the discrete Hilbert transform  
$(H_\mu 1) (z) =\int_{\T \setminus\{z\}} \frac{\,d\mu(\xi)}{1 -\bar\xi z}$ 
is bounded on $a(\mu)$: we have $|(H_\mu 1)(z)| \le C_\mu$ for all $z \in a(\mu)$. 
\end{itemize}
\end{Thm}
 The necessity of conditions $(a)$ and $(b)$ in Theorem \ref{t1} is well-known. I would like to thank A.~D.~Baranov who tell me the fact that condition $(c)$ is necessary as well. The proof of sufficiency part in Theorem \ref{t1} relies on a characterization of one-component inner functions in terms of their derivatives which is due to A.~B.~Aleksandrov~\cite{Al99}. 

\medskip 

\subsection{The main result} 
Having a description of the Clark measures of one-component inner functions, we now turn back to formula \eqref{eq36}. For a measure $\mu$ with properties~$(a)-(c)$ define the space $\bmo(\mu)$ by
$$
\bmo(\mu) = \left\{b \in L^1(\mu):\;  \|b\|_{\mu^*} = \sup_{\Delta}  \frac{1}{\mu(\Delta)}\int_{\Delta} | b - \langle b\rangle_{\Delta, \mu} | \, d\mu < \infty \right\},
$$ 
where $\Delta$ runs over all arcs of $\T$ with non-zero mass $\mu(\Delta)$ and $\langle b\rangle_{\Delta, \mu} = \frac{1}{\mu(\Delta)}\int_\Delta b\,d\mu$ is the standard integral mean of $b$ on $\Delta$. 
The following theorem is the main result of the paper.
\begin{Thm}\label{t2}
Let $\theta$ be a one-component inner function and let $\sigma_\alpha$ be its Clark measure. We have $(\Kl \cap z H^1)^* = \bmosa$. That is, for every continuous linear functional $\Phi$ on~$\Kl \cap zH^1$ there exists a function $b\in \bmosa$ such that $\Phi = \Phi_b$, where
\begin{equation}\label{eq24}
\Phi_b: F \mapsto \int_{\T} F b \,d\sigma_\alpha, \quad F \in \Kl \cap zH^\infty.
\end{equation}
Conversely, for every function $b \in \bmosa$ the functional~$\Phi_b$ is the densely defined continuous linear functional on $\Kl \cap zH^1$ with norm comparable to~$\|b\|_{\sa^*}$.
\end{Thm}
Every measure $\mu$ with properties $(a)$, $(b)$ from Theorem~\ref{t1} generates the doubling metric space $\bigl(\supp\mu,\,| \cdot |,\,\mu \bigr)$ in the sense of R.~Coifman and G.~Weiss~\cite{CoW83}. For such measures $\mu$ we have $H^{1}_{at}(\mu)^* = \bmo(\mu)$, where $H^{1}_{at}(\mu)$ is the corresponding atomic Hardy space,
\begin{equation}\label{eq25}
H^{1}_{at}(\mu) = \left\{\sum\nolimits_{k} \lambda_k a_k:\, a_k \mbox{ are $\mu$-atoms,} \;  \sum\nolimits_{k} |\lambda_k|< \infty \right\}.
\end{equation}
By a $\mu$-atom we mean a complex-valued function $a \in L^\infty(\mu)$ supported on an arc~$\Delta$ of~$\T$, with $\|a\|_{L^{\infty}(\mu)} \le 1/\mu(\Delta)$, and such that $\langle a\rangle_{\Delta, \mu} = 0$. The norm of $f \in H^{1}_{at}(\mu)$ is the infinum of $\sum_k|\lambda_k|$ over all possible representations $f = \sum\nolimits_{k} \lambda_k a_k$ of $f$ as a sum of $\mu$-atoms. We see from Theorem~\ref{t1} that Theorem~\ref{t2} admits the following equivalent reformulation. 
\addtocounter{Thm}{-1}
\renewcommand{\theThm}{\arabic{Thm}$'$} 
\begin{Thm}\label{t2prime}
Let $\mu$ be a measure with properties $(a) -(c)$. Then $f \in H^{1}_{at}(\mu)$ if and only if $f$ admits the analytic continuation to the unit open disk $\D$ as a function $F \in \Kl \cap zH^1$, where $\theta$ is the inner function with the Clark measure $\sigma_\alpha = \mu$. Moreover, such a function $F$ is unique and the norms $\|f\|_{H^{1}_{at}(\mu)}$, $\|F\|_{L^1(\T)}$ are comparable.
\end{Thm}
\renewcommand{\theThm}{\arabic{Thm}} 
For the counting measure $\mu = \delta_\Z$ on the set of integers $\Z$ Theorem \ref{t2prime} follows from the results by C.~Eoff \cite{Eo95}, S.~Boza and M.~Carro~\cite{BoC98}. They proved that $f \in H^{1}_{at}(\Z)$ if and only if $f$ admits the analytic continuation to the complex plane~$\C$ as a function from the \hbox{Paley-Wiener} space~$PW^1_{[0,2\pi]}$. It seems difficult to adapt the technique of \cite{BoC98} (where convolution operators were used to relate $H^{1}_{at}(\Z)$ and $\Re H^1(\R)$) for the general measures $\mu$ with properties $(a)-(c)$. Instead we give a complex-analytic proof based on the Cauchy-type formula
\begin{equation}\label{eq54}
\int_\Delta F(\xi)\,d\sa(\xi) = \oint_{\Gamma} \frac{F(z)/z}{1- \bar \alpha \theta(z)} \,dz,
\end{equation}
where $\Delta$ is an arc of $\T$, $\Gamma$ is a simple closed contour  in $\C$ which intersects $\T$ at the endpoints of~$\Delta$, and $F \in \Kl \cap z H^1$. Once we have a good estimate for the function $\frac{F(z)/z}{1- \bar \alpha \theta(z)}$ on $\Gamma$, formula \eqref{eq54} gives us an upper bound for the mean~$\langle F\rangle_{\Delta,\sa}$ on the arc~$\Delta$. 
Then we can use a standard Calder\'{o}n-Zigmund decomposition to obtain the representation of $F$ as a sum of atoms with respect to the measure~$\sa$.  The idea of using a contour integration is taken from the classical proof of atomic decomposition of $\Re(zH^1)$, where the contour $\Gamma$ comes from the Lusin-Privalov construction. In our situation we have to modify this construction so that the contour $\Gamma$ does not approach the subsets of the unit disk $\D$ where the function $|\alpha - \theta|$ is small.  

\medskip

\subsection{Truncated Hankel operators} 
One of important applications of the classical Fefferman duality theorem is the boundedness criterium for Hankel operators on the Hardy space $H^2$. Theorem \ref{t1} yields a similar criterium for truncations of Hankel  operators to coinvariant subspaces of $H^2$.  

\medskip

Let $\theta$ be an inner function and let $\Kth$ be the corresponding coinvariant subspace~\eqref{eq1} of the Hardy space $H^2$. Denote by $P_{\bar \theta}$ the orthogonal projection in $L^2(\T)$ to the subspace 
$\ov{z\Kth} = \{f \in L^2(\T):\, f= \ov{zg}, \; g \in \Kth\}$.  The truncated Hankel operator with symbol $\phi \in L^2(\T)$ is the densely defined operator $\Gamma_\phi: \Kth \to \ov{z\Kth}$,
\begin{equation}\label{eq75}
\Gamma_\phi: f \mapsto P_{\bar \theta} (\phi f), \quad f \in K_{\theta}^{\infty}.
\end{equation}
The symbol $\phi$ of $\Gamma_\phi$ is not unique. However, it is easy to check that every truncated Hankel operator on~$\Kth$ has the unique ``standard'' symbol $\phi \in \ov{\Kthtwo \cap z H^2}$, which plays the same role as the antianalytic symbol of a Hankel operator on $H^2$. 

\medskip

Two special cases of truncated Hankel operators are of traditional interest in the operator theory. If $\theta = z^n$, then the operators defined by \eqref{eq75} are classical Hankel matrices of size $n \times n$. Indeed, in this situation the space $\Kth$ consists of analytic polynomials of degree at most $n-1$ and  the entries of the matrix of $\Gamma_\phi$ in the standard bases of $\Kth$ and $\ov{z \Kth}$ depend only on the difference $k-l$: we have $(\Gamma_\phi z^k, \ov{z}^{l+1}) = \hat \phi(-k-l-1)$ for $0 \le k,l \le n-1$. Similarly, for the inner function $\theta_a: z \mapsto e^{iaz}$ in the upper half-plane $\C_+ = \{z \in \C: \;\Im z> 0\}$ the corresponding coinvariant subspace $K^{2}_{\theta_a}$ of the Hardy space $H^2(\C_+)$ can be identified with the  Paley-Wiener space 
$\pw^{2}_{[0, a]}$; truncated Hankel operators on $\pw^{2}_{[0, a]}$ are unitarily equivalent to the Wiener-Hopf convolution operators on the interval $[0, a]$, see \cite{Roch87, Carl11}. 

\medskip

The question for which symbols $\phi \in L^2(\T)$ the truncated Hankel operator $\Gamma_\phi$ is bounded on $\Kth$ (and how to estimate its operator norm in terms of $\phi$) admits several equivalent reformulations. It has been studied in \cite{Sar67, NikFarf03, Roch87, Carl11, BCFMT, BBK11}, see the discussion in Section \ref{s4}. Most of known results are Nehary-type theorems: under certain restrictions they affirm the existence of a bounded symbol for a bounded truncated Hankel/Toeplitz operator with control of the norms. Until now, the only $\bmo$-type criterium for truncated Hankel operators was known. In 2011, M.~Carlsson~\cite{Carl11} proved that a Hankel operator $\Gamma_\phi$ on $\pw^{2}_{[0, \pi]}$ with standard symbol~$\phi$ is bounded if and only if the sequence $\{\phi(n)\}_{n \in \Z}$ lies in the space $\bmo(\Z)$. Recall that we have $\pw^{2}_{[0, \pi]} = \mathcal K^{2}_{\theta_\pi}$ for the special one-component inner function $\theta_\pi: z \mapsto e^{i\pi z}$ in the upper half-plane~$\C_+$. The counting measure $\delta_\Z$ on $\Z$ can be regarded as the Clark measure $\nu_1$ for the inner function $\theta_\pi^2$ (for every inner function $\theta$ the Clark measures of $\theta^2$ will be denoted by $\nua$; from \eqref{eq20} we see that $\nua = (\sa + \sigma_{- \bar \alpha})/2$, $|\alpha| =1$).  Therefore the folowing result is a generalization of the criterium by M.~Carlsson.
\begin{Thm}\label{t3} Let $\theta$ be a one-component inner function, and let $\nu_\alpha$ be the Clark measure of the inner function $\theta^2$. The truncated Hankel operator $\Gamma_\phi: \Kth \to \ov{z\Kth}$ with standard symbol $\phi$ is bounded if and only if $\phi \in \bmo(\nua)$. Moreover, we have 
\begin{equation}\label{eq28}
c_1 \|\phi\|_{\nua^*} \le \|\Gamma_\phi\| \le c_2 \|\phi\|_{\nua^*},
\end{equation}
for some constants $c_1$, $c_2$ depending only on the inner function $\theta$. 
\end{Thm}
Similarly, one can describe compact truncated Hankel operators in terms of their standard symbols: we have $\Gamma_\phi \in S_\infty$ if and only if $\phi \in \vmo(\nua)$, see Section~\ref{s4}. 

\medskip

Theorem \ref{t3} for the inner function~$\theta = z^n$ yields the following interesting corollary for finite Hankel matrices.
\begin{Cor}\label{Cor1}
Let $\Gamma = (\gamma_{j+k})_{0 \le k, j \le n-1}$ be a Hankel matrix of size $n \times n$; consider its standard symbol $\phi = \gamma_0 \bar z + \gamma_1 \bar z^2 + \ldots \gamma_{2n-2} \bar z^{2n-1}$. We have 
\begin{equation}\label{eq2}
c_1 \|\phi\|_{\mu_{2n}^{*}} \le \|\Gamma\| \le c_2 \|\phi\|_{\mu_{2n}^{*}},
\end{equation}
where the constants $c_1, c_2$ do not depend on $n$ and  
$\mu_{2n} = \frac{1}{2n}\sum \delta_{\sqrt[2n]{1}}$ is the Haar measure on the group $\{\xi \in \C: \xi^{2n} = 1\}$. 
\end{Cor}
Corollary \ref{Cor1} implies the boundedness criterium for the standard Hankel operators on $H^2$. Recall that the Hankel operator $H_\phi: H^2 \to \ov{z H^2}$ with symbol $\phi \in L^2(\T)$ is densely defined by
$$H_\phi: f\mapsto P_-(\phi f), \quad f\in H^\infty,$$ 
where $P_-$ denotes the orthogonal projection in $L^2(\T)$ to $\ov{z H^2}$. It follows from the classical Fefferman duality theorem that $H_\phi$ is bounded if and only if its antianalytic symbol $P_- \phi$ lies in $\bmo(\T)$. Moreover, the operator norm of $H_\phi$ is comparable to $\|P_- \phi\|_{*}$, the norm  of $P_- \phi$ in $\bmo(\T)$. Taking the limit in \eqref{eq2} as $n \to \infty$ one can prove the estimate 
$c_1 \|\phi\|_{*} \le \|H_\phi\| \le c_2 \|\phi\|_{*}$ for every antianalytic polynomial $\phi$. This is already sufficient to obtain the general version of the boundedness criterium for Hankel operators on $H^2$, see details in Section \ref{s4}.
\section{Proof of Theorem \ref{t1}}\label{s2}
\subsection{Preliminaries} Given an inner function $\theta$, denote by $\rho(\theta)$ its boundary spectrum, that is, the set of points $\zeta \in \T$ such that $\liminf_{z \to \zeta, \, z \in \D} |\theta(z)| = 0$. In this paper we always assume that $\rho(\theta) \neq \T$, because this is so for one-component inner functions and for functions satisfying condition $(a)$ in Theorem \ref{t1} (see Lemma \ref{l4} below).  As is well-known, the function $\theta$ admits the analytic continuation from the open unit disk~$\D$ to the open domain $\D \cup G_\theta$, where $G_\theta = \bigl(\T\setminus\rho(\theta)\bigr)\cup \{z: \; |z| > 1,\; \theta(1/\bar z)\neq 0\}$. The analytic continuation is given by 
\begin{equation}\label{eq51}
\theta (z) = \frac{1}{\,\ov{\theta(1/\bar z)}\,}, \quad z \in G_\theta.
\end{equation}
Moreover, $\D \cup G_\theta$ is the maximal domain to which $\theta$ can be extended analytically. We need the following known lemma. 
\begin{Lem}\label{l4}
Let $\theta$ be an inner function with the Clark measure $\sa$, $|\alpha| = 1$. Then $\rho(\theta) = \rho(\sigma_\alpha)$. A point $z \in \T \setminus \rho(\theta)$ belongs to $\supp\sa$ if and only if $\theta(z) = \alpha$. Moreover, in the latter case we have $z \in a(\sa)$ and $\sa\{z\} = |\theta'(z)|^{-1}$.
\end{Lem}
\beginpf As is easy to see from formula \eqref{eq20}, we have
\begin{equation}\label{eq30}
\frac{\alpha + \theta(z)}{\alpha - \theta(z)} = \int_\T \frac{1 + \bar \xi z}{1 - \bar \xi z}\, d\sigma_\alpha(\xi) + i \Im \frac{\alpha + \theta(0)}{\alpha - \theta(0)}, \qquad z \in \D \cup G_\theta.
\end{equation}
Since $\theta$ is analytic on~$\D \cup G_\theta$, a point $z \in \T \setminus \rho(\theta)$ belongs to $\supp \sa$ if and only if $\theta(z) = \alpha$, and in the latter case there is no other points of $\supp\sa$ in a small neighbourhood of~$z$. Hence $z \in a(\sa)$ and we see from \eqref{eq30} that 
$$\sa\{z\} = (\bar \alpha z \theta'(z))^{-1} = |\theta'(z)|^{-1}.$$ It follows that $\T \setminus \rho(\theta) \subset \T \setminus \rho(\sigma_\alpha)$. For every $z \in \T \setminus \rho(\sigma_\alpha)$ either $z$ is an isolated atom of $\sigma_\alpha$ or $z \notin \supp\sigma_\alpha$. In both cases formula \eqref{eq30} shows that the function $\theta$ admits the analytic continuation from $\D$ to a small neighbourhood of~$z$. Hence $z \in \T \setminus \rho(\theta)$ and we have $\rho(\theta) = \rho(\sigma_\alpha)$. \qed

\medskip

The following result is in \cite{Al99}, see Theorem 1.11 and Remark 2 after its proof.
\begin{NMT}[A.~B.~Aleksandrov] An inner function $\theta$ is one-component if and only if it satisfies the following conditions:
\begin{itemize}
\item[(A1)] $m(\rho(\theta)) = 0$ and $|\theta'|$ is unbounded on every open arc $\Delta \subset \T\setminus \rho(\theta)$ such that $\ov{\Delta} \cap \rho(\theta) \neq \emptyset$; 
\item[(A2)] $\theta$ satisfies the estimate $|\theta''(\xi)| \le C |\theta'(\xi)|^2$ for all $\xi \in \T \setminus \rho(\theta)$.
\end{itemize} 
\end{NMT}

\subsection{Proof of Theorem \ref{t1}}  Essentially, we will show that conditions $(a) -(c)$ in Theorem \ref{t1} are equivalent to conditions $(A1)$, $(A2)$ above.

\medskip

\subsection*{Necessity} Let $\theta$ be a one-component inner function and let $\sa$ be its Clark measure. By Lemma \ref{l4} we have $\rho(\theta) = \rho(\sa)$. It was proved in \cite{Al99} that $m(\rho(\theta)) = 0$ and $\sa(\rho(\theta)) = 0$. Hence $\sa$ is a discrete measure with isolated atoms and we have  $m(\supp \sa) = 0$.  Let $\Delta$ be a connected component of the set $\T \setminus \rho(\sa) = \T \setminus \rho(\theta)$. By property $(A1)$ the argument of $\theta$ on $\Delta$ is a monotonic function  unbounded near both endpoints of~$\Delta$. It follows that the arc $\Delta$ contains infinitely many points $\xi_k$ such that $\theta(\xi_k) = \alpha$. Enumerate these points clockwise by integer numbers. We see from Lemma \ref{l4} that $\xi_k \in a(\sa)$ for all $k \in \Z$ and every atom $\xi_k$ has two neighbours $\xi_{k-1}$, $\xi_{k+1}$. This shows that the measure $\sa$ satisfies condition~$(a)$. 
The fact that $\sa$ satisfies condition $(b)$ follows from Lemma~5.1 of~\cite{BD11}.   
Now check condition $(c)$. Fix an atom~$\xi_0 \in a(\sa)$. 
From \eqref{eq30} we see that 
\begin{equation}\label{eq31}
\frac{1}{1 - \bar \alpha \theta(z)} = \int_\T \frac{d\sa(\xi)}{1 - \bar \xi z} + c_\alpha, \qquad z \in \D \cup G_\theta,
\end{equation}
where $c_\alpha = \alpha \ov{\theta(0)}/(1 - \alpha \ov{\theta(0)})$.
Hence,
\begin{equation}\notag
\begin{aligned}
(H_{\sa} 1)(\xi_0) + c_\alpha = \lim_{z \to \xi_0} \left(\frac{1}{1 - \bar \alpha \theta(z)} - \frac{\sa\{\xi_0\}}{1 - \bar \xi_0 z}\right).
\end{aligned}
\end{equation}
Consider the analytic function  $k_{\xi_0}: z \mapsto \frac{1 - \bar \alpha \theta (z)}{1 - \bar \xi_0 z}$ on the domain $\D \cup G_\theta$. We have 
\begin{equation}\label{eq45}
\begin{aligned}
(H_{\sa} 1)(\xi_0) + c_\alpha &= \lim_{z \to \xi_0} \frac{1}{1 - \bar \xi_0 z}\left(\frac{1}{k_{\xi_0}(z)}  - \frac{1}{k_{\xi_0}(\xi_0)}\right)\\
 &=-  \frac{\xi_0 k'_{\xi_0}(\xi_0)}{k_{\xi_0}^{2}(\xi_0)} = -\frac{\alpha \theta''(\xi_0)}{2\theta'(\xi_0)^2}.
\end{aligned}
\end{equation}
From here and the estimate in $(A2)$ we see that $H_{\sa} 1$ is bounded on $a(\sa)$. Surprisingly simple relation \eqref{eq45} between the discrete Hilbert transform $H_{\sa} 1$ and the inner function $\theta$ is the key observation in the proof.

\medskip

\subsection*{Sufficiency} Let $\mu$ be a measure with properties $(a) - (c)$. Construct the inner function~$\theta$ with the Clark measure $\sigma_\alpha = \mu$. To prove that $\theta$ is a one-component inner function we will check conditions $(A1)$ and $(A2)$.

\medskip

By Lemma \ref{l4} we have $\rho(\theta) = \rho(\sa)$. Hence $m(\rho(\theta)) = 0$ by property $(a)$ of the measure $\sa$. Let $\Delta$ be an open arc of $\T$ such that $\Delta \subset \T\setminus \rho(\theta)$ and $\bar \Delta \cap \rho(\theta) \neq \emptyset$. Then it follows from property~$(a)$ of the measure $\sa$ that $\Delta$ contains infinitely many atoms of $\sa$. Since $\sa$ is finite and $\sa\{\xi\} = |\theta'(\xi)|^{-1}$ for every $\xi\in a(\sa)$, the function $|\theta'|$ cannot be bounded on $\Delta$. This gives us condition~$(A1)$.

\medskip

Condition $(A2)$ is more delicate. To check it we need the following lemma. 

\begin{Lem}\label{l9} Assume that the Clark measure $\sa$ of an inner function $\theta$ has properties $(a)-(c)$. Then there exists a number $\kappa >0$ such that for every $\xi \in a(\sa)$ the set $D_\xi(\kappa) = \{z \in \C: |\xi - z| \le \kappa\sa\{\xi\}\}$ is contained in $\D \cup G_\theta$ and we have
\begin{equation}\label{eq46}
\frac{1}{2\sa\{\xi\}} \le \left|\frac{\alpha - \theta(z)}{\xi - z}\right| \le  \frac{2}{\sa\{\xi\}}
\end{equation}
for all $z \in D_\xi(\kappa)$.
\end{Lem}
\beginpf Pick an atom $\xi_0 \in a(\sa)$ and rewrite formula~\eqref{eq31} in the following form: 
\begin{equation}\label{eq33}
\frac{1}{1 - \bar \alpha \theta(z)} = \int_{\T \setminus \{\xi_0\}} \frac{d\sa(\xi)}{1 - \bar \xi z} + 
\frac{\sa\{\xi_0\}}{1 - \bar \xi_0 z} + c_\alpha, \qquad z \in \D\cup G_\theta.
\end{equation}
We have
\begin{equation}\notag
\left|\int_{\T \setminus \{\xi_0\}} \frac{d\sa(\xi)}{1 - \bar \xi z}\right| \le |(H_{\sa} 1) (\xi_0)| + \int_{\T \setminus \{\xi_0\}} \frac{|\xi_0 - z|\,d\sa(\xi)}{|\xi - z|\cdot|\xi - \xi_0 |}.
\end{equation}
By property $(c)$, $|(H_{\sa} 1) (\xi_0)| \le C_{\sa}$. Put $\kappa^* = (2B_{\sa})^{-1}$. For for $\xi \in a(\sa) \setminus\{\xi_0\}$ and $z \in D_{\xi_0}(\kappa^*)$ we have $|\xi_0 - z| \le \kappa^* \sa\{\xi_0\} \le |\xi_0 - \xi|/2$ by property $(b)$ of the measure $\sa$, which gives us the inequality $|\xi - z| \ge |\xi - \xi_0| - |\xi_0 - z| \ge \frac{1}{2}|\xi - \xi_0|$. It follows that for $z \in D_{\xi_0}(\kappa^*)$ we have 
\begin{equation} \notag
\int_{\T \setminus \{\xi_0\}} \frac{|z - \xi_0|\,d\sa(\xi)}{|\xi - z|\cdot|\xi - \xi_0 |} \le 2\kappa^* \sa\{\xi_0\} \int_{\T \setminus \{\xi_0\}} \frac{d\sa(\xi)}{|\xi - \xi_0|^2}.
\end{equation}
Denote by $\Delta$ the closed arc of $\T$ with endpoints $\xi_{0\pm} \in a(\sa)$. Using property $(b)$, we obtain the estimate 
\begin{equation} \label{eq39}
\begin{aligned}
\int_{\T \setminus \{\xi_0\}} \frac{d\sa(\xi)}{|\xi - \xi_0|^2} 
&\le \int_{\T \setminus \Delta} \frac{d\sa(\xi)}{|\xi - \xi_0|^2} + \frac{2}{A_{\sa}^{2} \sa\{\xi_0\}} \\
&\le 2\pi B_{\sa} \int_{\T \setminus \Delta} \frac{dm(\xi)}{|\xi - \xi_0|^2} + \frac{2}{A_{\sa}^{2} \sa\{\xi_0\}}\\
&\le \frac{C_1}{\sa\{\xi_0\}}, 
\end{aligned}
\end{equation}
where $C_1$ is a constant depending only on the measure $\sa$. We now see from \eqref{eq33} that 
\begin{equation}\label{eq59}
\frac{1}{1 - \bar \alpha \theta(z)} = \frac{\sa\{\xi_0\}}{1 - \bar \xi_0 z} + f_{\xi_0}(z), \quad z\in D_{\xi_0}(\kappa_1)\cap (\D \cup G_\theta),
\end{equation}
where the function $|f_{\xi_0}|$ is bounded by the constant $C_2 = 2 \kappa^* C_1  + C_{\sa} + |c_\alpha|$. Take a number $\kappa \le \kappa^*$ such that $C_2 \le (2\kappa)^{-1}$. We have $D_{\xi_0}(\kappa) \subset D_{\xi_0}(\kappa^*)$ and 
$$|f_{\xi_0}(z)| \le \frac{1}{2}\left|\frac{\sa\{\xi_0\}}{1 - \bar \xi_0 z}\right|$$
for all $z \in D_{\xi_0}(\kappa)\cap (\D \cup G_\theta)$. From here and \eqref{eq59} we get on $D_{\xi_0}(\kappa)\cap (\D \cup G_\theta)$ the estimate
$$
\frac{1}{2}\left|\frac{\sa\{\xi_0\}}{1 - \bar \xi_0 z}\right| \le \left|\frac{1}{1 - \bar \alpha \theta(z)}\right| \le 2 \left|\frac{\sa\{\xi_0\}}{1 - \bar \xi_0 z}\right|,
$$
which shows that $\theta$ admits the analytic continuation to a neighbourhood of $D_{\xi_0}(\kappa)$ (that is, $D_{\xi_0}(\kappa) \subset \D \cup G_\theta$) and proves formula \eqref{eq46} for points $z \in D_{\xi_0}(\kappa)$. Since our choice of the number $\kappa$ is uniform with respect to $\xi_0 \in a(\sa)$, the lemma is proved.  \qed

\medskip 

\noindent {\it Notation.} In what follows we write $E_1 \lesssim E_2$ (correspondingly, $E_1 \gtrsim E_2$) for two expressions $E_1, E_2$ to mean that there is a positive constant $c_\theta$ depending only on the inner function $\theta$ such that $E_1 \le c_\theta E_2$  (correspondingly, $c_\theta E_1 \ge  E_2$). We will write $E_1 \asymp E_2$ if $E_1\lesssim E_2$ and $E_1 \gtrsim E_2$.

\medskip

\noindent We are ready to complete the proof of Theorem \ref{t1}. Differentiating \eqref{eq31} we get
\begin{equation}\label{eq91}
\begin{aligned}
&\frac{\bar\alpha\theta'(z)}{(1 - \bar\alpha\theta(z))^2} = \int_{\T} \frac{\bar \xi d\sa(\xi)}{(1 - \bar \xi z)^2},  \\
&\frac{\bar\alpha\theta''(z)}{(1 - \bar\alpha\theta(z))^2} + \frac{2\bar\alpha^2 \theta'(z)^2}{(1 - \bar \alpha \theta(z))^3} = 2\int_\T \frac{\bar \xi ^2d\sa(\xi)}{(1 - \bar \xi z)^3}.
\end{aligned}
\end{equation}
Pick a point $\xi_0 \in a(\sa)$. Let $D_{\xi_0}(\kappa)$ be the set from Lemma~\ref{l9}. Denote by $\partial D_{\xi_0}(\kappa)$ the boundary of $D_{\xi_0}(\kappa)$. By formula \eqref{eq46}, $|\alpha - \theta(z)| \ge \kappa/2$ on ~$\partial D_{\xi_0}(\kappa)$. Arguing as in the Lemma \ref{l9}, from \eqref{eq91} we obtain the estimates
\begin{equation}\label{eq57}
\begin{aligned}
&|\theta'(z)| \lesssim \frac{\sa\{\xi_0\}}{|1 - \bar \xi_0 z|^2} + \int_{\T \setminus \{\xi_0\}} \frac{d\sa(\xi)}{|1 -  \bar \xi z|^2} \lesssim \frac{1}{\sa\{\xi_0\}},\\
&|\theta''(z)| \lesssim |\theta'(z)|^2 +  \frac{\sa\{\xi_0\}}{|1 -  \bar \xi_0 z|^3} + \int_{\T \setminus \{\xi_0\}} \frac{d\sa(\xi)}{|1 -  \bar \xi z|^3} \lesssim \frac{1}{\sa\{\xi_0\}^2}
\end{aligned}
\end{equation}
for all $z \in \partial D_{\xi_0}(\kappa)$. By the maximum principle we have $|\theta''(z)| \lesssim 1/\sa\{\xi_0\}^2$ for all points $z \in D_{\xi_0}(\kappa)$.
On the unit circle $\T$ we have 
\begin{equation}\notag
\frac{\bar \xi}{(1 - \bar \xi z)^2} =  \frac{- \bar z}{|1 - \bar \xi z|^2}.
\end{equation}
From here and formula \eqref{eq91} we get for $z \in D_{\xi_0}(\kappa) \cap \T$ the estimate  
\begin{equation}\label{eq34}
|\theta'(z)| = \sa\{\xi_0\} \left|\frac{1 - \bar\alpha\theta (z)}{1 - \bar \xi_0 z}\right|^2 + \int_{\T \setminus \{\xi_0\}} \left|\frac{1 -\bar \alpha \theta(z)}{1 -  \bar \xi z}\right|^2\,d\sa(\xi) \gtrsim \frac{1}{\sa\{\xi_0\}}.
\end{equation}
Combining \eqref{eq57} and \eqref{eq34} we see that $|\theta''/\theta'^2| \lesssim 1$ on $D_{\xi_0}(\kappa) \cap \T$. It remains to obtain the same estimate for points $z \in \T \setminus\rho(\theta)$ that do not belong to the union of the sets $D_{\xi}(\kappa)$, $\xi \in a(\sa)$, from Lemma \ref{l9}. Take such a point $z_0$.
We claim that $|\theta(z_0) - \alpha| \ge \kappa/2$. Indeed, assume the converse and find the connected component $\Delta$ of the set $\{\zeta \in \T\setminus\rho(\theta): |\theta(\zeta) - \alpha| < \kappa/2\}$ containing the point $z_0$. Since the argument of the inner function $\theta$ is monotonic on $\Delta$, there exists a point $\xi \in \Delta$ such that $\theta(\xi) = \alpha$. By Lemma~\ref{l4} we have $\xi \in a(\sa)$. Next, from \eqref{eq46} we see that $|\alpha - \theta(z)| \ge \kappa/2$ for both points in $\T\cap \partial D_{\xi}(\kappa)$. Since $\Delta$ is connected this yields the inclusion $\Delta \subset D_\xi (\kappa)$ which gives us the contradiction with $z_0 \notin D_\xi (\kappa)$. Thus, we proved the inequality $|\theta(z_0) - \alpha| \ge \kappa/2$. Let $\xi_{z_0}$ be the nearest point to $z_0$ in $a(\sa)$. We have $\kappa \sa\{\xi_{z_0}\} \le |\xi_{z_0} - z| \le B_{\sa} \sa\{\xi_{z_0}\}$. These two estimates imply~\eqref{eq57} and~\eqref{eq34} for $z = z_0$ and $\xi_0 = \xi_{z_0}$.
It follows that $|\theta''(z_0)/\theta'(z_0)^2| \lesssim 1$ and $\theta$ satisfies condition~$(A2)$. \qed

\medskip 

\noindent{\bf Remark.} Lemma 5.1 in \cite{BD11} and formula \eqref{eq45} show that for every one-component inner function $\theta$ there exist positive constants $A_\theta$, $B_\theta$, $C_\theta$ such that 
$A_\theta \le A_{\sa}$, $B_{\sa} \le B_\theta$, $C_{\sa} \le C_\theta$ for all  Clark measures~$\sa$ of $\theta$. Also, it follows from Lemma~5.1 in~\cite{BD11} that $|\theta'(z)| \asymp 1/\sa\{\xi\}$ for all $\xi \in a(\sa)$ and all $z \in \T$ between the neighbours~$\xi_\pm$ of $\xi$ in $a(\sa)$. In particular, we have $\sigma_\beta(\Delta) \asymp \sa(\Delta) \asymp m(\Delta)$ for all $\beta$ with $|\beta| = 1$ and all arcs $\Delta$ of the unit circle $\T$ containing at least two atoms of the measure $\sa$.

\section{Proofs of Theorem \ref{t2} and Theorem \ref{t2prime}}\label{s3}
\noindent We first prove Theorem \ref{t2prime}. 
The following result is classical, for the proof see \cite{Cl72} or Chapter 9 in \cite{Cima2006}.
\begin{NMT}(D.\ N.\ Clark)\label{l5}
Let $\theta$ be an inner function and let $\sigma_\alpha$ be its Clark measure. The natural embedding $V_\alpha:  \Kth \into L^2(\sa)$ defined on the reproducing kernels of the space $\Kth$ by  
$V_\alpha\left(\frac{1 - \ov{\theta(\lambda)} \theta}{1 - \bar \lambda z}\right) = \frac{1 - \ov{\theta(\lambda)} \alpha}{1 - \bar \lambda z}$
can be extended to the whole space~$\Kth$ as the unitary operator from $\Kth$ to $L^2(\sa)$. For every $f \in L^2(\sigma_\alpha)$ the function 
\begin{equation}\label{eq40}
F(z) = \int_\T f(\xi) \frac{1 - \bar \alpha \theta(z)}{1 - \bar \xi z}\,d\sigma_\alpha(\xi)
\end{equation}
in the unit disk $\D$ belongs to $\Kth$ and $V_\alpha F = f$ as elements in~$L^2(\sigma_\alpha)$. 
\end{NMT}
It worth be mentioned that A.\ G.\ Poltoratski \cite{Polt93} established the existence of angular boundary values $\sa$-almost everywhere on $\T$ for all functions in the space~$\Kth$, thus proving that the unitary operator $V_\alpha$ in Clark theorem acts as the natural embedding on the whole space $\Kth$. In our situation this follows from a very simple argument, see Lemma \ref{l131} in Section \ref{s32}. 

The embedding $V_\alpha: \Kthp \into L^p(\sa)$ defined on the linear span of the reproducing kernels of~$\Kthp$ might be unbounded for $1 \le p < 2$ and might have the unbounded inverse $V_{\alpha}^{-1}: L^p(\sa) \into \Kthp$ for $2 < p \le \infty$, see Section 3 in \cite{Al89}. However, the situation is ideal for the one-component inner functions~$\theta$, as following results show:
\begin{itemize}
\item $V_\alpha \Kthp \subset L^p(\sa)$ for $1 < p < \infty$ -- A. L. Volberg, S. R. Treil \cite{VT86};
\item $V_\alpha \Kthp = L^p(\sa)$ for $1 < p < \infty$ -- A. B. Aleksandrov \cite{Al89}; 
\item $V_\alpha \Kthp \subset L^p(\sa)$ for $0 < p \le 1$ -- A. B. Aleksandrov \cite{Al99}.
\end{itemize}
Theorem \ref{t2prime} says that $V_\alpha \Kl = H^{1}_{at}(\sa)$ for every one-component inner function $\theta$. We are ready to prove its easy part -- the inclusion $V_\alpha \Kl \supset H^{1}_{at}(\sa)$.

\medskip

\subsection{Proof of the part ``$\Rightarrow$'' in Theorem \ref{t2prime}} Let $\mu$ be a measure on the unit circle~$\T$ with properties $(a)-(c)$. Take a complex number $\alpha$ of unit modulus and construct the one-component inner function $\theta$ with the Clark measure $\sa = \mu$. We want to show that every function $f \in H^1_{at}(\sigma_\alpha)$ admits the analytic continuation to the open unit disk $\D$ as a function $F \in \Kl \cap z H^1$ with $\|F\|_{L^1(\T)} \lesssim \|f\|_{H^{1}_{at}(\sa)}$. At first, assume that $f$ is a $\sa$-atom supported on an arc $\Delta \subset \T$ with center~$\xi_c$. Then $f \in L^2(\sigma_\alpha)$ and the function $F$ in formula \eqref{eq40} lies in the space $\Kth \subset \Kl$ by Clark theorem. Since $\int_{\T} f \,d\sa = 0$, we  have $F(0) = 0$. Moreover, we see from Lemma \ref{l4} that $F(\xi) = f(\xi)$ for all $\xi \in a(\sa)$. Let us check that the norm of $F$ in $L^1(\T)$ is bounded by a constant depending only on the inner function $\theta$. By Aleksandrov desintegration theorem (see \cite{AB87}  or Section 9.4 in \cite{Cima2006}), we have 
\begin{equation}
\int_\T |F|\,dm = \int_{\T}\int_\T|V_\beta F(\xi)|\,d\sigma_\beta(\xi)\,dm(\beta).
\end{equation}
Fix a complex number $\beta \neq \alpha$ of unit modulus. We claim that $\|V_\beta F\|_{L^1(\sigma_\beta)} \lesssim 1$. Denote by $2\Delta$ the arc of $\T$ with  center $\xi_c$ such that $m(2 \Delta) = 2m(\Delta)$ (in the case where $m(\Delta) \ge 1/2$ put $2\Delta = \T$). Break the integral $\int_\T |V_\beta F|\,d\sigma_\beta$ into two parts,  
\begin{equation}\label{eq41}
\int_\T|V_\beta F(\xi)|\,d\sigma_\beta(\xi) = \int_{2\Delta}|V_\beta F(\xi)|\,d\sigma_\beta(\xi) + \int_{\T\setminus 2\Delta}|V_\beta F(\xi)|\,d\sigma_\beta(\xi).
\end{equation}
By Clark theorem we have $\|V_\beta F\|_{L^2(\sigma_\beta)} = \|F\|_{L^2(\T)} = \|V_\alpha F\|_{L^2(\sigma_\alpha)}$. Moreover, we have $\|V_\alpha F\|_{L^2(\sigma_\alpha)} \le 1/\sqrt{\sigma_\alpha(\Delta)}$ because the function $V_\alpha F = f$ is a $\sa$-atom supported on the  arc $\Delta$. This yields the inequality
\begin{equation}\label{eq83}
\int_{2\Delta}|V_\beta F(\xi)|\,d\sigma_\beta(\xi) \le \sqrt{\sigma_\beta(2\Delta)} \cdot \|V_\beta F\|_{L^2(\sigma_\beta)}  \le \sqrt{\sigma_\beta(2\Delta)/\sigma_\alpha(\Delta)}.
\end{equation}
Note that the arc $\Delta$ contains at least two points in $a(\sa)$ because $f$ has zero $\sa$-mean on $\Delta$. Hence $\sa(\Delta) \asymp m(\Delta)$ and $\sigma_\beta(2\Delta) \asymp m(2\Delta)$, see remark after the proof of Theorem~\ref{t1}. This shows that $\int_{2\Delta}|V_\beta F(\xi)|\,d\sigma_\beta(\xi) \lesssim 1$. Let us now estimate the second term in \eqref{eq41}. Take a point $z \in a(\sigma_\beta)\setminus 2\Delta$. Using the fact that $f$ is a $\sa$-atom we obtain the estimate
\begin{equation}\label{eq42}
\begin{aligned}
|V_\beta F(z)|
&=\left|\int_\Delta f(\xi) \frac{1 - \bar \alpha \beta}{1 - \bar \xi z}\,d\sigma_\alpha(\xi)\right|\\
&=\left|\int_\Delta f(\xi)\left(\frac{1 - \bar \alpha \beta}{1 - \bar \xi z} - \frac{1 - \bar \alpha \beta}{1 - \bar \xi_c z}\right)\,d\sigma_\alpha(\xi)\right|\\
&\le 2\int_\Delta |f(\xi)|\left|\frac{\xi - \xi_c}{(1 - \bar \xi z)(1 - \bar \xi_c z)}\right| \,d\sigma_\alpha(\xi)\\
&\le \frac{2\pi m(\Delta)}{|z - \xi_c|^2} \cdot \sup_{\xi \in \Delta}\left|\frac{z - \xi_c}{z - \xi}\right| \cdot \int_\Delta |f(\xi)|\,d\sigma_\alpha(\xi)\\ &\le \frac{4\pi m(\Delta)}{|z - \xi_c|^2}.\\
\end{aligned}
\end{equation}
From here we get
\begin{equation}\label{eq422}
\int_{\T\setminus 2\Delta} |V_\beta F(z)|
\lesssim m(2\Delta) \cdot \int_{\T\setminus 2\Delta}\frac{d\sigma_\beta(z)}{|z - \xi_c|^2} \lesssim 1.
\end{equation}
Hence the norm of $F$ in $L^1(\T)$ is bounded by a constant depending only on $\theta$. Now take an arbitrary function $f \in H^{1}_{at}(\sa)$ and consider its representation $f = \sum \lambda_k f_k$, where $f_k$ are $\sa$-atoms and $\sum_k |\lambda_k| \le 2 \|f\|_{H^{1}_{at}(\sa)}$. Let $F_k$ be the functions in $\Kth$ such that $V_\alpha F_k = f_k$. Then the sum $\sum \lambda_k F_k$ converges absolutely in $L^1(\T)$ to a function $F \in \Kl$ and we have $\|F\|_{L^1(\T)} \lesssim  \|f\|_{H^{1}_{at}(\sa)}$. From formula \eqref{eq40} we get
$$F(z) =  \sum_k \lambda_k \int_{\T} f_k(\xi) \frac{1 - \bar \alpha \theta(z)}{1 - \bar \xi z}\,d\sa(\xi) = \int_\T f(\xi) \frac{1 - \bar \alpha \theta(z)}{1 - \bar \xi z}\,d\sa(\xi), \quad z \in \D.$$
Since $f \in L^1(\T)$, this formula determines the analytic continuation of $F$ to the domain $\D \cup G_\theta$. From Lemma~\ref{l4} we see that $F(\xi) = f(\xi)$ for all $\xi \in a(\sa)$. \qed

\subsection{Preliminaries for the proof of the part ``$\Leftarrow$'' in Theorem \ref{t2prime}}\label{s32}
Let $\theta$ be a one-component inner function with the Clark measure $\sa$. Introduce positive constants $\tilde A_{\sa}$, $\tilde B_{\sa}$ such that
$$\tilde A_{\sa} m[\xi, \xi_\pm] \le  \sa\{\xi\} \le \tilde B_{\sa}  m[\xi, \xi_\pm], \quad \xi \in a(\sa).$$
Here $[\xi, \xi_-]$, $[\xi, \xi_+]$ are the closed arcs of $\T$ with endpoints $\xi, \xi_\pm \in a(\sa)$ such that the corresponding open arcs $(\xi, \xi_\pm)$ do not intersect $\supp\sa$. Take a positive number $\kappa \le (2 \tilde B_{\sa})^{-1}$ for which estimate \eqref{eq46} holds true. Denote by $D_{\sa}(\kappa)$ the union   of the sets $D_{\xi}(\kappa)$, $\xi \in a(\sa)$ from Lemma \ref{l9}.
\begin{Lem}\label{l17}
For every arc $\Delta$ of $\T$ containing at least one atom of the measure~$\sa$ we have $m(\Delta) \le (2/\tilde A_{\sa}) \sa(\Delta)$. If $\Delta$ contains two or more atoms of $\sa$, we have  $\sa(\Delta) \le 4\tilde B_{\sa} m(\Delta \setminus D_{\sa}(\kappa))$. In particular, the sets $D_{\xi}(\kappa)$, $\xi \in a(\sa)$ are disjoint. 
\end{Lem}
\beginpf It is sufficient to prove the statement in the case where $\Delta$ contains only finite number of atoms of $\sa$. Enumerate the atoms clockwise: $\xi_1, \dots, \xi_n$. Find the neighbours of $\xi_1, \xi_n$ in $a(\sa) \setminus \Delta$ and denote them by $\xi_{0}$ and $\xi_{n+1}$, correspondingly. We have 
$$
m(\Delta) \le \sum_{k=0}^{n} m[\xi_k, \xi_{k+1}] \le \frac{2}{\tilde  A_{\sa}} \sum_{k=1}^{n} \sa\{\xi_k\} = \frac{2}{\tilde A_{\sa}}\sa(\Delta).
$$
In the case where $n \ge 2$ we have
\begin{equation}\notag
\sa(\Delta) = \sum_{k=1}^{n} \sa\{\xi_k\} \le 2\tilde B_{\sa} m(\Delta) 
\le 2\tilde B_{\sa} m(\Delta \setminus D_{\sa}(\kappa)) + \tilde B_{\sa} \kappa \sa(\Delta).
\end{equation} 
Now use the assumption $\kappa \le(2\tilde B_{\sa})^{-1}$ and  get $\sa(\Delta) \le 4\tilde B_{\sa}m(\Delta \setminus D_{\sa}(\kappa))$. \qed

\begin{Lem}\label{l7}
There exists $\eps> 0$ such that $|\alpha - \theta(z)| \ge \eps$ for all $z \in \D\setminus D_{\sa}(\kappa)$. 
\end{Lem}
\beginpf Let $\delta \in (0,1)$ be a number such that the set $\Omega_\delta = \{z \in \D: |\theta(z)|< 1\}$ is 
connected. The set $$\Omega_{\delta,\frac{1}{\delta}} = \{z \in \D\cup G_\theta : \delta < |\theta(z)|< 1/\delta\}$$
is at most countable union of the open connected components, $\mathcal O_k$. It was proved by B.~Cohn \cite{Co82} that the restriction of  the inner function  $\theta$ to each of the sets $\mathcal O_k$ is a covering map from $\mathcal O_k$ to the ring $R_\delta = \{z \in \C: \delta < |z| < 1/\delta\}$. Take a positive number  $\eps< \min(\kappa/2, 1 -\delta)$. We claim that every connected component $E$ of the set $L_\eps = \{z \in \D \cup G_\theta: \; |\alpha - \theta(z)| < \eps\}$ contains an atom of $\sa$. Indeed, we have $E \subset \mathcal O_k$ for some index $k$ because $L_\eps \subset \Omega_{\delta,\frac{1}{\delta}}$. Since $\theta$ is a covering map from $\mathcal O_k$ to $R_\delta$, there exists a number $\eps_1$ (which can be taken to be less than $\eps$) such that the preimage of $\{\zeta:\; |\alpha - \zeta| < \eps_1\}$ under $\theta$ on $\mathcal O_k$ is at most countable union of the open disjoint sets $\mathcal O_{km} \subset \mathcal O_k$ and $\theta$ is a homeomorphism from $\mathcal O_{km}$ to $\{\zeta:\; |\zeta - \alpha| < \eps_1\}$ for every $m$. 
 By the minimum principle, $\inf_{z \in E}|\theta(z) - \alpha| = 0$. It follows that $E \cap \mathcal O_{km} \neq \emptyset$ for some index $m$. Since $E$ is connected and $|\theta - \alpha| < \eps_1 < \eps$ on $\mathcal O_{km}$, we have $\mathcal O_{km} \subset E$. But every set $\mathcal O_{km}$ contains the unique point $\xi$  with $\theta(\xi) = \alpha$. By Lemma \ref{l4}, $\xi \in a(\sa)$ and thus $E \cap a(\sa) \neq \emptyset$. To prove the lemma it is sufficient to show that $E \subset D_\xi(\kappa)$. For every $z \in E \cap D_\xi(\kappa)$ we get from~\eqref{eq46} the estimate 
$$
|\xi - z| \le 2|\alpha - \theta(z)|\sa\{\xi\} \le 2\eps \sa\{\xi\}. 
$$
Hence $E$ does not intersect the circle $\{z \in \C: |z - \xi| = r\sa\{\xi\}\}$ for 
every $r \in (2\eps, \kappa)$. Since the set $E$ is connected this yields the desired inclusion $E \subset D_\xi(\kappa)$. \qed 

\medskip

\begin{Lem}\label{l131}
Let $\theta$ be an inner function with $\rho(\theta) \neq \T$. Then every function in~$\Kl$ admits the analytic continuation from the unit disk $\D$ to the domain $\D \cup G_\theta$. Consequently, if $\sa(\rho(\theta)) = 0$ for a Clark measure $\sa$ of $\theta$, then every function in~$\Kl$ has a trace on the set $a(\sa) \subset \D \cup G_\theta$ of full measure $\sa$.     
\end{Lem}
\beginpf For every function $F \in \Kl$ we have $\bar \theta F \in \ov{z H^1}$ on $\T$. Hence, 
\begin{equation}\label{eq60}
F(z) = \int_\T F(\xi)\frac{1 - \theta(z)\ov{\theta(\xi)}}{1 - z \bar\xi}\,dm(\xi), \qquad z \in \D.
\end{equation}
Extend the inner function $\theta$ to the domain $\D \cup G_\theta$ by formula \eqref{eq51}. The right hand side of \eqref{eq60} then determines the analytic continuation of the function $F$ to $\D \cup G_\theta$. By Lemma \ref{l4} we have $\rho(\sa) = \rho(\theta)$ which completes the proof. \qed

\medskip

\begin{Lem}\label{l6}
Let $\theta$ be an inner function and let $G \in \Kl \cap zH^1$. Then there exist  functions $G_1, G_2 \in \Kl \cap z H^1$ such that $G = G_1 + i G_2$ and $G_{1,2} = \theta \ov{G}_{1,2}$ on $\T \setminus \rho(\theta)$. Moreover, we have $\|G_{1,2}\|_{L^1(\T)} \le \|G\|_{L^1(\T)}$. 
\end{Lem}
\beginpf Consider the function $\tilde G = \theta \ov{G}$ on the unit circle $\T$. We have 
$$\tilde G \in \theta(\ov{z H^1} \cap z \bar \theta H^1) = 
\bar z \theta \ov{H^1} \cap z H^1 = \Kl \cap zH^1.$$ 
This shows that $G$ can be continued to the open unit disk $\D$ as a function from the space $\Kl \cap z H^1$. Now put $G_1 = (G + \tilde G)/2$, $G_2 = (G - \tilde G_1)/2i$ and obtain the desired representation. \qed

\medskip

\subsection{Proof of the part ``$\Leftarrow$'' in Theorem \ref{t2prime}.} Let $\mu$ be a measure on the unit circle~$\T$ with properties $(a)-(c)$ and let $|\alpha| = 1$. Consider the one-component inner function~$\theta$ with the Clark measure $\sa = \mu$. Take a function $F \in \Kl \cap zH^1$.  By Lemma~\ref{l131}, $F$ is analytic on the domain $\D\cup G_\theta$. Denote by $f$ its  trace on the set $a(\sa) \subset \D \cup G_\theta$ of full measure $\sa$. Our aim is to prove that $f \in H^{1}_{at}(\sa)$ and $\|f\|_{H^{1}_{at}(\sa)} \lesssim \|F\|_{L^1(\T)}$. At first, assume that $F \in \Kth \cap z H^2$ and $F = \theta \ov{F}$ on $\T \setminus \rho(\theta)$. We will need the following modification of the Lusin-Privalov construction (see Section~III.D in \cite{Koos98} for the standard one). Consider the non-tangential maximal function of $F$,
$$F^*(\xi) = \sup_{z \in \Lambda_\xi} |F(z)|, \quad \xi \in \T,$$
where $\Lambda_\xi$ denotes the convex hull of the set $\{\xi\}\cup \{z\in\D: |z| \le 1/\sqrt{2}\}$. Put 
$$S_{F}(\lambda) = \ov{\D}\setminus\{z \in \ov{\D}: z \in \Lambda_\xi \mbox{ for some } \xi\in \T \mbox{ with }  F^*(\xi) < \lambda\}.$$ 
Let $D_{\sa}(\kappa)$ be the set defined at the beginning of Section \ref{s32}. By Lemma \ref{l7}, we have $|\alpha - \theta| \ge \eps$ on $\D \setminus D_{\sa}(\kappa)$. Denote by $R_F(\lambda)$ the union of those connected components of the set $S_F(\lambda) \cup D_{\sa}(\kappa)$ for which we have $E \cap S_{F}(\lambda) \neq \emptyset$ and $E \cap D_{\sa}(\kappa) \neq \emptyset$. The sets $R_F(\lambda)$ are closed and have the following properties:
\begin{itemize}
\item[(1)] If $\lambda_1 < \lambda_2$, then $R_F(\lambda_2) \subset R_F(\lambda_1)$;
\item[(2)] $|F(z)| \le \lambda$ for $\sa$-almost all points $z \in \T\setminus R_F(\lambda)$;
\item[(3)] $|F(z)| \le \lambda$ and $|\alpha - \theta(z)| \ge \eps$ for $z \in \partial R_F(\lambda) \cap \D$.
\end{itemize}
More special properties of the sets $R_F(\lambda)$ are collected in the following lemma.
\begin{Lem}\label{l8}
Let $E$ be a connected component of the set $R_F(\lambda)$. Put $\gamma = \partial E \cap \D$ and $\Delta = \partial E \cap \T$. There exist constants $c_4$, $c_5$, $c_6$ depending only on $\theta$ such that
\begin{itemize}
\item[(4)] $\gamma$ is a rectifiable curve with length $|\gamma| \le c_4 \sa(\Delta)$;
\item[(5)] $\sa(\Delta) \le c_5m(\Delta \cap S_F(\lambda))$ if $E$ contains  at least two atoms of $\sa$;
\item[(6)] $\frac{1}{\sa(\Delta)}\left|\int_{\Delta}f\,d\sigma_\alpha\right| \le c_6 \lambda$.
\end{itemize}
One can take $c_4 = 40 /\tilde A_{\sa}$, $c_5 = 4\tilde B_{\sa}$, $c_6 = 60/(\eps \tilde A_{\sa})$.
\end{Lem}
\beginpf By the construction and Lemma \ref{l17} we have
$$|\gamma| \le (\sqrt{2}+ \pi/2)|\Delta| \le 20 m(\Delta) \le 40 \sa(\Delta)/\tilde A_{\sa}.$$
In the case where the arc $\Delta$ contains at least two atoms of the measure $\sa$ Lemma~\ref{l17} gives us the estimate 
$$\sa(\Delta) \le 4\tilde B_{\sa} m(\Delta \setminus D_{\sa}(\kappa)) \le   4\tilde B_{\sa} m(\Delta \cap S_F(\lambda)).$$
Les us check property~$(6)$. At first, assume that $\gamma \cap \rho(\theta) = \emptyset$. Then we have $\gamma \cap \supp\sa = \emptyset$ by the construction. For $z \in \C$ with $|z|\ge 1$ denote $z^* = 1/\bar z$ and put $\gamma^* = \{z \in \C: z^* \in \gamma\}$. The set $\Gamma = \gamma \cup \gamma^*$ is a rectifiable curve in $\C$ with length $|\Gamma| \le 3|\gamma|$. Let us check that 
\begin{equation}\label{eq61}
\left|\frac{F(z)/z}{1- \bar \alpha \theta(z)}\right| \le 2\eps^{-1}\lambda, \quad z \in \Gamma \cap (\D \cup G_\theta).
\end{equation}
For $z \in \gamma$ we have $|z|\ge 1/\sqrt{2}$, $|F| \le \lambda$, $|\alpha - \theta| \ge \eps$ and therefore \eqref{eq61} holds. The function $z \mapsto \ov{F(z^*)/\theta(z^*)}$ is analytic on the interior of $G_\theta$ and coincides with the function $F$ on $G_\theta \cap \T = \T\setminus \rho(\theta)$ (recall that $F$ admits the analytic continuation to the domain $\D\cup G_\theta$ by Lemma \ref{l131} and $F = \theta \ov{F}$ on $\T \setminus \rho(\theta)$ by the assumption). By the the uniqueness of the analytic continuation  we have $F(z) = \ov{F(z^*)/\theta(z^*)}$ for all $z \in G_\theta$. Now take a point $z \in G_\theta$ and compute 
$$
\frac{F(z)/z}{1- \bar \alpha \theta(z)} = \frac{\ov{z^* F(z^*)/\theta(z^*)}}{1- \ov{\alpha/\theta(z^*)}} = \ov{\frac{z^* F(z^*)}{\theta(z^*) - \alpha}}.
$$
This yields estimate \eqref{eq61} for $z \in \gamma^* \cap G_\theta$. Next, we claim that
\begin{equation}\label{eq53}
\int_\Delta f(\xi)\,d\sa(\xi) = - \frac{1}{2\pi i}\oint_{\Gamma} \frac{F(z)/z}{1- \bar \alpha \theta(z)} \,dz.
\end{equation}
Indeed, using formula \eqref{eq40} for the function $F/z \in \Kth$ we obtain
\begin{equation}\label{eq62}
\begin{aligned}
\oint_{\Gamma} \frac{F(z)/z}{1- \bar \alpha \theta(z)}\,dz
&=
\oint_{\Gamma} \frac{1}{1- \bar \alpha \theta(z)} \int_{\T} \bar \xi f(\xi) \frac{1- \bar \alpha \theta(z)}{1 - \bar \xi z} \,dz\,d\sa(\xi)\\
&= \int_\T f(\xi) \oint_{\Gamma} \frac{1}{\xi - z} \,dz\,d\sa(\xi)= - 2\pi i \int_\T f(\xi) \chi_{\Delta}(\xi)\,d\sa(\xi),
\end{aligned}
\end{equation}
where $\chi_{\Delta}$ denotes the indicator of the set $\Delta$.
Note that change of the order of integration is possible because $\Gamma \cap \supp \sa = \emptyset$ and therefore all integrals in \eqref{eq62} are absolutely convergent. We now see from \eqref{eq61} and \eqref{eq53} that 
\begin{equation}\notag
\left|\int_\Delta f(\xi)\,d\sa(\xi)\right| \le 
(\pi\eps)^{-1} \lambda |\Gamma| 
\le 3(\pi\eps)^{-1} \lambda |\gamma| \le 3(\pi\eps)^{-1} c_3  \cdot  \lambda \cdot \sa(\Delta).
\end{equation}
This gives us property $(5)$ in the case where $\gamma \cap \rho(\theta) = \emptyset$. The general case can be reduced to just considered one by a small perturbation of the contour $\gamma$; use the fact that $f \in L^2(\sa)$ by Clark theorem and property $(a)$ of the measure $\sa$ from from Theorem \ref{t1}. \qed

\bigskip

Lemma \ref{l8} is the key argument in the proof of Theorem \ref{t2prime}. The rest of the proof is a standard Calder\'{o}n-Zigmund decomposition. We will follow the exposition in Section VII.E of \cite{Koos98}. For each $\lambda > 0$ the set $\Delta_F(\lambda) = R_F(\lambda) \cap \T$ is a union of closed disjoint arcs 
$\Delta_F^k(\lambda)$, $\Delta_F(\lambda) = \cup_{k \in I_\lambda}\Delta_F^k(\lambda)$. Consider the functions 
\begin{equation}\notag
G_\lambda = \left\{
\begin{aligned}
&f, \hspace{1.1cm}\xi \in \T\setminus \Delta_F(\lambda),\\
&\langle f \rangle_{\Delta_F^k(\lambda),\sa}, 
\;\xi \in \Delta_F^k(\lambda),\\
\end{aligned}
\right. 
\quad 
B_\lambda = \left\{
\begin{aligned}
&0, \hspace{1.8cm} \xi \in \T\setminus \Delta_F(\lambda),\\
&f - \langle f \rangle_{\Delta_F^k(\lambda),\sa}, 
\; \xi \in \Delta_F^k(\lambda).\\
\end{aligned}
\right.
\end{equation}
By Lemma \ref{l8} we have $|G_\lambda| \le c_{6}\lambda$ $\sa$-almost everywhere on $\T$. The function $B_\lambda$ has zero $\sa$-mean on each arc $\Delta_{F}^{k}(\lambda)$, $k \in I_\lambda$. For every integer $n \in \Z$ set $g_n = G_{2^n}$ and $b_n = B_{2^n}$. Fix a number $N_0 \in \Z$ such that 
$$2^{N_0} < \inf_{|z| \le 1/\sqrt{2}}|F(z)| 
\le 2^{N_0+1}.$$
Note that $\Delta_F(2^{N_0}) = \T$. By formula \eqref{eq40},  
$$
g_{N_0} = \frac{1}{\sa(\T)} \int_\T f \, d\sa = \frac{F(0)}{\sa(\T) (1 - \bar \alpha \theta(0))} = 0. 
$$
Since $f$ is finite at each point $\xi \in a(\sa)$ we have $f(\xi) =  g_N(\xi)$ for every sufficiently big number $N$. Hence 
\begin{equation}\label{eq85}
f(\xi) = \sum_{n = N_0}^{\infty}\bigl(g_{n+1}(\xi) - g_{n}(\xi)\bigr), \quad \xi \in a(\sa),
\end{equation}
where the sum converges pointwise (in fact,  only finite number of summands in \eqref{eq85} are non-zero for every $\xi \in a(\sa)$). Note that $f = b_n + g_n$ and  $g_{n+1} - g_{n} = b_{n} - b_{n+1}$ for all $n \ge N_0$.
Let $I'_{2^n}$ be the set of indexes $k \in I_{2^n}$ such that the set $\Delta_{F}^{k}(2^n)$ contains at least two atoms of the measure~$\sa$. 
The function $g_{n+1} - g_{n}$ vanishes \hbox{$\sa$-almost} everywhere on each of the sets $\Delta_{F}^{k}(2^n)$, $ k \in I_{2^n} \setminus I'_{2^n}$. Indeed, for such index $k$ we have by the construction. Hence 
$g_n(\xi) = g_{n+1}(\xi) = f(\xi)$ because the $\sa$-mean of $f$ on any arc containing the only point $\xi \in a(\sa)$ equals~$f(\xi)$. Define 
\begin{equation}\notag
\tilde{a}_{n,k} = \chi_{\Delta_{F}^{k}(2^n)} (b_n - b_{n+1}), \qquad n \ge N_0, \quad k \in I'_{2^n},
\end{equation}
where $\chi_{\Delta_{F}^{k}(2^n)}$ is the indicator of the set $\Delta_{F}^{k}(2^n)$. The functions $\tilde{a}_{n,k}$ have zero \hbox{$\sa$-mean} on $\T$. Indeed, let $I$ denote the set of indexes $m$ such that $\Delta_F^m(2^{n+1}) \subset \Delta_F^k(2^{n})$ (note that $\Delta_F(2^{n+1}) \subset \Delta_F(2^{n})$ by property $(1)$ of the sets $R_F(\lambda)$). Then 
$$
\int_\T \tilde a_{n,k}\,d\sa = \int_{\Delta_{F}^k(2^n)} (b_{n} - b_{n+1})\,d\sa = -\sum_{m \in I} \int_{\Delta_F^m(2^{n+1})} b_{n+1}\,d\sa = 0.
$$
Also, we have $|\tilde a_{n,k}| \le |g_n| + |g_{n+1}| \le 3 c_6 \cdot  2^{n}$ on $\T$ for every $n \ge N_0$ and $k \in I'_{2^n}$. Now put
$$
a_{n,k} = \frac{\tilde a_{n,k}}{3 c_6 \cdot 2^n \cdot  \sa(\Delta_{F}^{k}(2^n))}, \qquad n\ge N_0, \quad k \in I'_{2^n},
$$
and observe that $a_{n,k}$ are atoms with respect to the measure $\sa$. It follows from formula \eqref{eq85} that
\begin{equation}\label{eq69}
f(\xi) = \sum_{n \ge N_0} \sum_{k \in I'_{2^n}} \lambda_{n,k} a_{n,k}(\xi), \qquad \xi \in a(\sa), 
\end{equation}
where $\lambda_{n,k} = 3 c_6 \cdot 2^n \cdot  \sa(\Delta_{F}^{k}(2^n))$ and the sum is convergent 
pointwise. Since the set $a(\sa)$ has full measure $\sa$ it remains to check that 
\begin{equation}\label{eq68}
\sum_{n \ge N_0} \sum_{k \in I'_{2^n}} \lambda_{n,k} \lesssim \|F\|_{L^1(\T)}. 
\end{equation}
By Lemma \ref{l8} we have 
$$\sa(\Delta_{F}^{k}(2^n)) \le c_5 m(\Delta_{F}^{k}(2^n) \cap S_{F}(2^n))$$
for every $n \ge N_0$ and $k \in I'_{2^n}$. Hence, 
\begin{equation}\notag
\begin{aligned}
&\sum_{n \ge N_0} \sum_{k \in I'_{2^n}} 2^{n} \sa(\Delta_{F}^{k}(2^n)) 
\le c_5\sum_{n \ge N_0} \sum_{k \in I'_{2^n}} 2^{n} m(\Delta_{F}^{k}(2^n) \cap S_{F}(2^n)) \\
&\le c_5 \sum_{n \ge N_0} 2^{n} m(S_{F}(2^n) \cap \T) 
= c_5 \sum_{n \ge N_0} 2^{n} m(\{\xi \in \T: F^*(\xi) \ge 2^n\}).
\end{aligned}
\end{equation}
The last sum does not exceed
\begin{equation}\notag
\begin{aligned}
\sum_{n \ge N_0}& \sum_{l \ge 0} 2^{n} m\bigl(\{\xi \in \T:  2^{n+l} \le F^*(\xi) < 2^{n+l+1}\}\bigr) \le \\
&\le\sum_{l \ge 0} m\bigl(\{\xi \in \T: 2^{N_0+l}\le F^*(\xi) < 2^{N_0+l+1}\}\bigr) \sum\nolimits_{k = N_0}^{l} 2^{N_0+k}  \\  
& \le \sum_{l \ge 0} 2^{N_0 + l+ 1} \cdot m\bigl(\{\xi \in \T: 2^{N_0+l}\le F^*(\xi) < 2^{N_0+l+1}\}\bigr)\\
&\le 2\|F^*\|_{L^1(\T)} \le 2M\|F\|_{L^1(\T)},
\end{aligned}
\end{equation}
where $M$ denotes the norm of the maximal operator $F \mapsto F^*$ on $H^1$. Thus, inequality \eqref{eq68} holds with the constant $6c_5 c_6M$ and formula \eqref{eq69} gives us the atomic decomposition of the trace $f$ provided $F \in \Kth \cap zH^2$ and $F = \theta \bar F$. Now consider arbitrary function $F \in \Kl \cap zH^1$ with the trace $f$ on the set $a(\sa)$. Since $\Kth \cap zH^2$ is the dense subset of $\Kl \cap zH^1$ in norm of $L^1(\T)$ one can find functions $F_k \in \Kth \cap zH^2$ such that $F = \sum_k F_k$ and $\|F\|_{L^1(\T)} \ge \frac{1}{2} \sum_{k} \|F_k\|_{L^1(\T)}$. Let $G_{1,k}$, $G_{2,k}$ be the functions from Lemma \ref{l6} for $G = F_k$ and let $g_{1,k}$, $g_{2,k}$ be their traces on $a(\sa)$. We have $f(\xi) = \sum g_{1,k}(\xi) + i\sum g_{2,k}(\xi)$ for every $\xi \in a(\sa)$, see formula~\eqref{eq60}. It follows from the first part of the proof that $f$ admits the atomic decomposition with respect to the measure $\sa$ and we have $\|f\|_{H^{1}_{at}(\sa)} \le 24 c_5 c_6 M \|F\|_{L^1(\T)}$. \qed

\medskip

\subsection{Proof of Theorem \ref{t2}} Since $(\supp\sa,\, |\cdot|,\, \sa)$ is the doubling metric space, Theorem~\ref{t2prime} and Theorem B in \cite{CoW83} imply Theorem \ref{t2}. To make the paper more self-contained, we give a proof of this implication.  

\beginpf Let $\theta$ be a one-component inner function. We first remark that the integral in formula~\eqref{eq24} is correctly defined for $F \in K^{\infty}_{\theta}$ and $b \in \bmo(\sa)$. Indeed, by Lemma \ref{l4} and Lemma \ref{l131} every function $F \in \Kl$ has the trace $f$ on the set $a(\sa)$ of full measure $\sa$. If $F \in \Kl \cap zH^\infty$, then $f \in L^\infty(\sa)$. Since $\bmo(\sa) \subset L^1(\T)$ the integral in formula \eqref{eq24} converges absolutely. 

Consider a continuous linear functional $\Phi$ on $\Kl \cap zH^1$. Since $\Kth \subset \Kl$ and $\Kth \cap zH^2$ is the Hilbert space there exists a function $G \in \Kth \cap z H^2$ such that $\Phi(F) = \int_{\T} F\ov{G} \, dm$ for all $F \in \Kth \cap zH^2$. Denote by $b$ the restriction of the function $\ov{G}$ to the set $a(\sa)$ of full measure $\sa$. By Clark theorem, we have  $b \in L^2(\sa)$. Let us prove that $b  \in \bmosa$. For every function $F \in \Kth \cap zH^2$ we have
\begin{equation}\label{eq73}
\Phi(F) = \int_{\T} F\ov{G} \, dm = \int_{\T} F b \, d\sa = \Phi_b(F),    
\end{equation}
where use again Clark theorem. Take an arc $\Delta$ of $\T$ and consider the function $a_{0} \in L^\infty(\sa)$  such that $|a_0| = 1$, $a_0(b - \langle b \rangle_{\Delta, \sa}) = |a_0(b - \langle b \rangle_{\Delta, \sa})|$ $\sa$-almost everywhere on $\Delta$ and $a = 0$ $\sa$-everywhere off $\Delta$. Denote by $\chi_\Delta$ the indicator of the set $\Delta$. The function 
$$
a = \frac{1}{2\sa(\Delta)}(a_0 - \langle a_0 \rangle_{\Delta, \sa})\chi_\Delta
$$ 
is an atom with respect to the measure $\sa$ and we have
\begin{equation}\label{eq80}
\int_\T a b \,d\sa = \int_\Delta a (b - \langle b \rangle_{\Delta, \sa}) \,d\sa = \frac{1}{2\sa(\Delta)}\int_\Delta |b - \langle b \rangle_{\Delta, \sa}| \,d\sa.  
\end{equation}
By Theorem \ref{t2prime} the function $a$ can be continued analytically to $\D$ as a function $F_{a} \in \Kl \cap zH^1$ with $\|F_a\|_{L^1(\T)} \lesssim 1$. Since $a \in L^2(\sa)$, we have $F_a \in \Kth \cap zH^2$ by Clark theorem. Now it follows from \eqref{eq73} and \eqref{eq80} that 
$\|b\|_{\sa^*} \lesssim \|\Phi_b\|$. 

\medskip

Conversely, take a function $b \in \bmosa$ and consider the functional $\Phi_b$ densely defined on $\Kl \cap zH^1$ by formula~\eqref{eq24}. For every $\sa$-atom $a$ supported on an arc $\Delta$ we have 
\begin{equation}\label{eq84}
\left|\int_\T a b \,d\sa\right| = \left|\int_\Delta a (b - \langle b \rangle_{\Delta, \sa}) \,d\sa \right| \le \frac{1}{\sa(\Delta)}\int_\Delta |b - \langle b \rangle_{\Delta, \sa}| \,d\sa.  
\end{equation}
This shows that the functional $f \mapsto \int_{\T} f b \,d\sa$ is continuous on $H^{1}_{at}(\sa)$. By Theorem \ref{t2prime}, the restriction of every function $F \in \Kl \cap z H^1$ to $a(\sa)$ belongs to $H^{1}_{at}(\sa)$ and $\|F\|_{H^{1}_{at}(\sa)} \lesssim \|F\|_{L^1(\T)}$. Hence the functional $\Phi_b$ is continuous on $\Kl \cap zH^1$ and we see from \eqref{eq84} that $\|\Phi_{b}\| \lesssim \|b\|_{\sa^*}$. \qed

\section{Truncated Hankel and Toeplitz operators}\label{s4}
Let $\theta$ be an inner function. Denote by $P_\theta$ the orthogonal projection in $L^2(\T)$ to the subspace $\Kth$. The truncated Toeplitz operator $A_\psi: \Kth \to \Kth$ with symbol $\psi \in L^2(\T)$ is densely defined by
$$
A_\psi: f \mapsto P_\theta (\psi f), \quad f \in K^{\infty}_{\theta}.
$$
Truncated Toeplitz and Hankel operators are closely related. Indeed, the antilinear isometry $g \mapsto \bar z \theta \bar g$ on $L^2(\T)$ preserves the subspace $\Kth$ and for every $f, g \in K^{\infty}_{\theta}$ we have
\begin{equation}\label{eq78}
(A_\psi f, g) = (\psi f, g) = (\Gamma_{\bar \theta \psi} f,\ov{zg_1}),
\qquad g_1 = \bar z \theta \bar g.
\end{equation}
This shows that the operators $A_\psi$, $\Gamma_{\bar \theta \psi}$ are bounded (compact, of trace class, etc.) or not simultaneously and $\|A_\psi\| = \|\Gamma_{\bar\theta\psi}\|$. Below we briefly discuss some results  related to the boundedness problem for truncated Toeplitz operators. 

\medskip

We will say that the truncated Toeplitz operator $A_\psi$ has a bounded symbol~$\psi_1$ if $A_\psi = A_{\psi_1}$ for a function $\psi_1 \in L^\infty(\T)$. It can be shown all symbols of the zero truncated Toeplitz operator on $\Kth$ have the form $\ov{\theta g_1} + \theta g_2$, where $g_1, g_2 \in H^2$, see \cite{Sar07}. Hence the operator $A_\psi: \Kth \to \Kth$ has a bounded symbol if and only if the set $\psi + \ov{\theta H^2} + \theta H^2$ contains a bounded function on $\T$. 
Clearly, every truncated Toeplitz operator with bounded symbol is bounded. The following question arises: does every bounded truncated Toeplitz operator have a bounded symbol? 

\medskip

\subsection{Analytic symbols} In 1967, D.~Sarason \cite{Sar67} described the commutant $\{\Sth\}'$ of the restricted shift operator $\Sth: f \mapsto P_\theta (zf)$ on $\Kth$. He proved that a bounded operator~$A$ on $\Kth$ commutes with $\Sth$ if and only if there exists a function $\psi \in H^\infty$ such that $A = A_\psi$. Moreover, we have $\|A_\psi\| = \dist_{H^\infty}(\psi, \theta H^\infty)$ and one can choose the function $\psi$ so that $\|A\| = \|\psi\|_{H^\infty}$. This well-known theorem yields a boundedness criterium for truncated Toeplitz operators with analytic symbols. Indeed, for every $\psi \in H^2$ and $f \in K^{\infty}_{\theta}$ we have $A_\psi \Sth f = \Sth A_\psi f$. Hence the operator $A_\psi$ is bounded if and only if $A_\psi \in \{\Sth\}'$  which is equivalent to the existence of a function $\psi_1 \in H^\infty$ such that $A_\psi = A_{\psi_1}$ (in other words, we have $\psi + \theta h \in H^\infty$ for some $h \in H^2$). The equality $\|A_\psi\| = \dist_{H^\infty}(\psi, \theta H^\infty)$ for $\psi \in H^\infty$ leads to a short proof for the Nevanlinna-Pick interpolation theorem and its generalization, see \cite{Sar67}. 

It was observed by N.~K.~Nikolskii that many problems for truncated Toeplitz operators with analytic symbols can be easily reduced to the problems for usual Hankel operators on $H^2$. The reduction is based on the fact that for every $\psi \in H^2$ the operator $\bar \theta A_\psi P_{\theta}$ from $H^2$ to $\ov{z H^2}$ coincides with the Hankel operator $H_{\bar \theta \phi}$. In particular, the operator $A_\psi$ is bounded (compact, of trace class, etc.) if and only if so is the  operator $H_{\bar \theta \psi}$. Since Hankel operators on $H^2$ are well studied this observation immediately yields consequences for truncated Toeplitz operators. As an example, the operator $A_\psi$ on $\Kth$ with symbol $\psi \in H^2$ is compact if and only if $\bar \theta \psi \in C(\T) + H^2$, where $C(\T)$ denotes the algebra of continuous functions on the unit circle $\T$. For more information see Lecture 8 in \cite{Nik86} and Section 1.2 in \cite{PeBook}.

\medskip

\subsection{General symbols}
Until recently, a little was known about truncated Toeplitz operators with general symbols in $L^2(\T)$. For such operators the boundedness problem is more complicated. 

\medskip

In 1987, R.~Rochberg \cite{Roch87} proved that every bounded Toeplitz operator on the Paley-Wiener space $\pw^{2}_{[-a, a]}$ has a bounded symbol. Using the Fourier transform, he reduced the general case of the problem to consideration of the Toeplitz operators on $\pw^{2}_{[0, a]}$ with analytic symbols. Recently, M.~Carlsson \cite{Carl11} use a result from \cite{Roch87} to prove the boundedness criterium for Topelitz and Hankel operators on $\pw^{2}_{[-a, a]}$ in terms of $\bmo(\frac{\pi}{a}\Z)$, see Section \ref{s1}.

\medskip

Every finite Toeplitz matrix $A$ clearly have bounded symbols. However, the question concerning the best possible constant $c_{A}$ in the inequality 
$$\inf\{\|\psi\|_{L^\infty(\T)}: A_\psi = A\} \le c_{A} \cdot \|A\|$$ 
is nontrivial.
In 2001,  M.~Bakonyi and D.~Timotin proved that $c_A \le 2$ for every self-adjoint finite Toeplitz matrix $A$. As a corollary, we have $c_A \le 4$ for a general finite Toeplix matrix~$A$ that was improved to $c_A \le 3$ by L.~N.~Nikolskaya and Yu.~B.~Farforovskaya~\cite{NikFarf03} in 2003. Next, in 2007 D.~Sarason \cite{Sar07} compute $c_A = \pi/2$ for 
$A = \left(\begin{smallmatrix}0 & i \\ -i & 0\end{smallmatrix}\right)$ and proved that $c_A \le \pi/2$ for every $2\times 2$ self-adjoint Toeplitz matrix~$A$. In paper \cite{Vol04} A.~L.~Volberg discuss several approaches to the dual version of the problem of determining $\sup_{A} c_A$ over all finite Toepliz matrices $A$, which can be formulated in terms of weak factorizations of analytic polynomials.

\medskip
 
In 2010, A.~D.~Baranov, I.~Chalendar, E.~Fricain, J.~Mashreghi, and D.~Ti\-motin~\cite{BCFMT} constructed an inner function $\theta$ and a bounded truncated Toeplitz operator $A$ on~$\Kth$ that has no bounded symbols. Shortly after that in \cite{BBK11} appeared a description of coinvariant subspaces $\Kth$ on which every bounded truncated Toeplitz operator has a bounded symbol. The proof in \cite{BBK11} is based on a duality relation between the space of all bounded truncated Toeplitz operators on $\Kth$ and a special function space.  With help of formula \eqref{eq78} it is easy to reformulate the results of \cite{BBK11} for truncated Hankel operators. We do this below as a preparation to the proof of Theorem \ref{t3}.

\subsection{Duality for truncated Hankel operators}\label{s43} Let $\theta$ be an inner function. Consider the linear space 
$$
\Yth = \left\{\sum_{k = 0}^{\infty}x_k y_k, \; x_k \in \Kth, \,
y_k \in z\Kth, \;\;
\sum_{k = 0}^{\infty} \|x_k\|_{L^2(\T)}\|y_k\|_{L^2(\T)}
< \infty\right\}.
$$
As is easy to see, we have $\Yth \subset \Kths \cap zH^1$. Define the norm in $\Yth$ by
$$
\|h\|_{\Yth} =  \inf\left\{\sum_{k = 0}^{\infty} 
\|x_k\|_{L^2(\T)}\|y_k\|_{L^2(\T)}: \;
h = \sum_{k = 0}^{\infty}x_k y_k, \; x_k \in \Kth, \,
y_k \in z\Kth 
\right\}.
$$ 
With this norm $\Yth$ is a Banach space. Denote by $\Hth$ the linear space of all bounded truncated Hankel operators acting from $\Kth$ to $\ov{z \Kth}$. It follows from Theorem 4.2 in \cite{Sar07} that $\Hth$ is closed in the weak operator topology. Hence $\Hth$ is the Banach space under the standard operator norm and moreover it has a predual space. It follows from Theorem 2.3 of \cite{BBK11} that $\Yth^* = \Hth$. That is, for every continuous linear functional~$\Psi$ on $\Yth$ there exists the unique operator $\Gamma \in \Hth$ such that~$\Psi = \Psi_{\Gamma}$, where
\begin{equation}\label{eq76}
\Psi_\Gamma: h \mapsto \sum_{k = 0}^{\infty}(\Gamma x_k, \ov{y_k}), \qquad h \in \Yth, \quad h = \sum_{k = 0}^{\infty}x_k y_k.
\end{equation}
Conversely, for every operator $\Gamma \in \Hth$ the mapping $\Psi_{\Gamma}$ is the correctly defined continuous linear functional on the space $\Yth$ and we have $\|\Psi_\Gamma\| = \|\Gamma\|$. 

\medskip

With help of the equality $\Yth^* = \Hth$ the boundedness problem for truncated Hankel operators can be reformulated in terms of function theory. Indeed, now it is easy to see from Hahn-Banach theorem that every bounded truncated Hankel operator on $\Kth$ has a bounded symbol if and only if $\Yth$ is a closed subspace of $L^1(\T)$, in which case $\Yth$ coincides with $\Kths \cap zH^1$ as a set, see details in \cite{BBK11}. Note that if $\Yth = \Kths \cap zH^1$ as set, then the the norms $\| \cdot \|_{\Yth}$ and $\| \cdot \|_{L^1(\T)}$ are equivalent on $\Yth$. It was proved in \cite{BBK11} that $\Yth = \Kths \cap zH^1$ for every one-component inner function $\theta$.

\medskip

Thus, we see from the results of \cite{BBK11} and Theorem \ref{t2} that for every one-component inner function $\theta$ we have
$$
\Yth^* = \Hth, \quad \Yth = \Kths \cap z H^1, \quad (\Kths \cap z H^1)^* = \bmo(\nua),
$$
where $\nua$ is the Clark measure of the inner function $\theta^2$. It remains to combine this relations to obtain Theorem \ref{t3}.

\subsection{Proof of Theorem \ref{t3}} Let $\theta$ be a one-component inner function and let $\Gamma_\phi$ be a truncated Hankel operator on $\Kth$ with standard symbol $\phi \in \ov{\Kthtwo \cap zH^2}$; we do not assume now that the operator $\Gamma_\phi$ is bounded. For every function $h = \sum_{k = 0}^{\infty}x_k y_k$ in $\Yth \cap L^\infty(\T)$ we have 
\begin{equation}\label{eq79}
\Psi_{\Gamma_\phi}(h) = \sum_{k = 0}^{\infty}(\Gamma_\phi x_k, \ov{y_k}) = 
\int_\T \phi \sum_{k = 0}^{\infty}x_k y_k \,dm = \int_\T h \phi\,dm = \int_\T h \phi\,d\nu_\alpha,
\end{equation}
where the last equality follows from Clark theorem for the inner function $\theta^2$. 
We see that $\Psi_{\Gamma_\phi}$ coincides on $\Yth \cap L^\infty(\T)$ with the functional 
$$\Phi_{\phi}: h \mapsto \int_\T h \phi\,d\nu_\alpha, \quad h \in K^{1}_{\theta^2} \cap z H^\infty.$$
Since the inner function $\theta^2$ is one-component the Banach spaces $\Yth$ and \hbox{$\Kths \cap z H^1$} coincide as sets and their norms are equivalent. It follows that the densely defined functionals $\Psi_{\Gamma_\phi}: \Yth \to \C$ and $\Phi_{\phi}: \Kths \cap z H^1 \to \C$ are bounded or not  simultaneously and $\|\Psi_{\Gamma_\phi}\| \asymp \|\Phi_{\phi}\|$, where the constants involved depend only on~$\theta$. By Theorem \ref{t2} for the inner function $\theta^2$ the functional $\Phi_{\phi}$ is bounded if and only if $\phi \in \bmo(\nua)$, and in the latter case we have $\|\Phi_{\phi}\| \asymp \|\phi\|_{\nua^*}$. Now result follows from the equality $\|\Gamma_\phi\| = \|\Psi_{\Gamma_\phi}\|$. \qed

\subsection{Compact truncated Hankel operators} Let $\mu$ be a measure on $\T$ with properties $(a)-(c)$. For every $b \in \bmo(\mu)$ define
$$
M_\eps(b) = \sup\left\{\frac{1}{\mu(\Delta)} \int_\Delta |b - \langle b \rangle_{\Delta, \mu}| \,d\mu, \; \Delta \mbox{ is an arc of $\T$ with } 0< \mu(\Delta) \le \eps \right\}.
$$ 
Consider the space $\vmo(\mu) = \{b \in \bmo(\mu): \lim_{\eps \to 0} M_\eps(b) =0\}$ of functions of vanishing mean oscillation with respect to the measure $\mu$. It can be shown that $\vmo(\mu)$ is the closure in $\bmo(\mu)$ of the set of all finitely supported sequences.
\begin{Prop}\label{p1}
Let $\theta$ be a one-component inner function, and let $\nu_\alpha$ be the Clark measure of the inner function $\theta^2$. The truncated Hankel operator $\Gamma_\phi: \Kth \to \ov{z\Kth}$ with standard symbol $\phi$ is compact if and only if $\phi \in \vmo(\nua)$.
\end{Prop}
\beginpf It follows from Theorem 2.3 of \cite{BBK11} that $(\Hth \cap S_\infty)^* = \Yth$, where $S_\infty$ denotes the ideal of all compact operators acting from $\Kth$ to $\ov{z \Kth}$. Hence a bounded truncated Hankel operator $\Gamma$ on $\Kth$ is compact if and only if the functional $\Psi_{\Gamma}$ in~\eqref{eq76} is continuous in the weak$^*$ topology on $\Yth$. Let $\phi$ be the standard symbol of the operator $\Gamma_\phi$. By Corollary 2.5 in \cite{BBK11} and Theorem \ref{t2prime} we have 
$$\Yth = \Kths \cap z H^1, \qquad V_\alpha(\Kths \cap z H^1) = H^{1}_{at}(\nua).$$ 
From formula \eqref{eq79} we see that the operator $\Gamma_\phi$ is compact if and only if the restriction of $\phi$ to $a(\nua)$ generates the weak$^*$ continuous functional $\Phi_\phi: f \mapsto \int f \phi \,d\nua$ on the space $H^{1}_{at}(\nua)$.  
For any doubling measure $\mu$ we have $\vmo(\mu)^* = H^{1}_{at}(\mu)$, see Theorem 4.1 in \cite{CoW83}. It follows that $\Gamma_{\phi} \in S_\infty$ if and only if $\phi \in \vmo(\nua)$.\qed

\subsection{Functions in $\Kth$ of bounded mean oscillation} Theorem~\ref{t3} provides the following description of functions in $\Kth \cap \bmo(\T)$. 
\begin{Prop}\label{p2}
Let $\theta$ be a one-component inner function and let $\phi \in \Kth$. Then we have $\phi \in \Kth \cap \bmo(\T)$ if and only if $\phi \in \bmo(\nua)$, where $\nua$ is the Clark measure of the inner function $\theta^2$.
\end{Prop}
\beginpf A function $\phi \in H^2$ belongs to the space $\bmo(\T)$ if and only if the Hankel operator $H_{\bar \phi}: H^2 \to \ov{z H^2}$ is bounded, see Theorem 1.2 in Chapter 1 of \cite{PeBook}. Assume that $\phi \in \Kth$ and consider the truncated Hankel operator $\Gamma_\phi: \Kth \to \ov{z\Kth}$. For every function $f \in K^\infty_\theta$ we have $\bar \phi f \in \bar\theta H^2$. Hence, 
$$H_{\bar \phi} f = P_{-} (\bar \phi f) = P_{\bar \theta} (\bar \phi f) = \Gamma_{\bar \phi} f, \quad f \in K^{\infty}_{\theta}.$$
Also, $H_\phi f = 0$ for all $f \in \theta H^\infty$. Therefore the operators $H_{\bar \phi}$ and $\Gamma_{\bar \phi}$ are bounded or not simultaneously and $\|H_{\bar \phi}\| = \|\Gamma_{\bar \phi}\|$. Now the result follows from Theorem~\ref{t3}. \qed

\medskip

\subsection{Finite Hankel and Toeplitz matrices}
\noindent Let $\Gamma = (\gamma_{j+k})_{0 \le j,k \le n-1}$ be a Hankel matrix of size $n \times n$. Associate with $\Gamma$ the antianalytic polynomial 
$$\phi = \gamma_0 \bar z + \gamma_1 \bar z^2 + \ldots \gamma_{2n-2} \bar z^{2n-1}.$$
For the inner function $\thn =  z^{n}$ the space $\Kthn$ consists of analytic polynomials of degree at most $n-1$. Consider the truncated Hankel operator $\Gamma_\phi: \Kthn \to \ov{z\Kthn}$,
$$\Gamma_\phi: f \mapsto P_{\bar{\theta}_n} (\phi f), \quad f \in \Kthn.$$
We have $(\Gamma_\phi z^j, {\bar z}^{k+1}) = \gamma_{j+k}$ for every $0 \le j, k \le n-1$. It follows that 
the matrix~$\Gamma$ as the operator on $\C^{n}$ is unitary equivalent to the operator $\Gamma_\phi$. Analogously, the Toeplitz matrix $A = (\alpha_{j-k})_{0 \le j,k \le n-1}$ is unitarily equivalent to the truncated Toeplitz operator $A_\psi: \Kthn \to \Kthn$ with symbol 
$$\psi = \alpha_{-(n-1)}\bar z^{n-1} + \ldots + \alpha_{n-1}z^{n-1}.$$
If moreover $\alpha_m = \gamma_{(n-1)-m}$ for all $m \in \Z$ with $|m|\le n-1$, then $\phi = \bar\theta_n \psi$ and  we have $\|\Gamma\| = \|A\|$ by formula \eqref{eq78}. Consider the measure
\begin{equation}\notag
\mu_{2n} = \frac{1}{2n}\sum \delta_{\sqrt[2n]{1}}   
\end{equation}
equally distributed at the roots of identity of order $2n$: $\supp \mu = \{\xi \in \T: \xi^{2n} = 1\}$. 
Let $c_{1, n}$,  $c_{2,n}$ be the best possible constants in the inequality   
\begin{equation}\label{eq81}
c_{1, n}\|\phi\|_{\mu_{2n}^{*}} \le \|\Gamma_\phi\| \le c_{2,n} \|\phi\|_{\mu_{2n}^{*}},
\end{equation}
where $\Gamma_\phi$ runs over all truncated Hankel operators on $\Kthn$, $\phi$ is the standard symbol of $\Gamma_\phi$. Corollary \ref{Cor1} of Theorem~\ref{t3} claims that the sequences $\{c_{1,n}^{-1}\}_{n \ge 1}$ and $\{c_{2,n}\}_{n \ge 1}$ are bounded. We prove this below.

\medskip

\noindent {\bf Proof of Corollary \ref{Cor1}.} We may assume that $n \ge 2$. It follows from Lemma~\ref{l4} that $\mu_{2n}$ is the Clark measure $\nu_1$ of the inner function $\thn^2 = z^{2n}$. This allows us to estimate the constants in formula \eqref{eq81} using the proofs of Theorem \ref{t2} and Theorem~\ref{t3}. Denote
\begin{equation}
\begin{aligned}
&d'_{n} = \sup\{\|h\|_{L^1(\T)},\; h \in \Kthns \cap z H^1, \; \|h\|_{H^{1}_{at}(\mu_{2n})} \le 1\};\\
&d''_{n} = \sup\{\|h\|_{Y_{\thn}},\; h \in \Kthns \cap z H^1, \; \|h\|_{L^1(\T)} \le 1\}.
\end{aligned}
\end{equation}
Let $\Gamma_\phi: \Kthn \to \ov{z \Kthn}$ be a truncated Hankel operator with standard symbol $\phi$. Consider the functional $\Psi: h \mapsto \int_{\T} h \phi \,d\mu_{2n}$ on the Banach space $Y_{\thn}$. From formula~\eqref{eq80} and the equality $\|\Psi\| = \|\Gamma_\phi\|$ (see Section \ref{s43}) we obtain
\begin{equation}
\begin{aligned}
\|\phi\|_{\mu_{2n}^{*}} &\le 2\sup\{|\Psi(h)|, \; \|h\|_{H^{1}_{at}(\mu_{2n})} \le 1\} \le 
2d'_{n}\sup\{|\Psi(h)|, \; \|h\|_{L^1(\T)} \le 1\} \\
&\le 2d'_{n} d''_{n}\sup\left\{|\Psi(h)|, \; \|h\|_{Y_{\thn}} \le 1 \right\} = 2d'_{n} d''_{n}\|\Gamma_\phi\|. 
\end{aligned}
\end{equation}
Hence, $c_{1, n}^{-1} \le 2 d'_{n} \cdot d''_{n}$. It follows from the results of Nikolskaya and Far\-fo\-rovs\-kaya \cite{NikFarf03} that $d''_{n} \le 3$, see also Section 1.2 in \cite{Vol04}. To estimate the constant $d'_{n}$ assume that the restriction of $f \in \Kthns \cap z H^1$ to $a(\mu_{2n})$ is a $\mu_{2n}$-atom supported on a closed arc $\Delta$ of the unit circle $\T$ with center~$\xi_c$ and endpoints in $a(\mu_{2n})$. Let $\{\nu_{\beta}^{n}\}_{|\beta| = 1}$ be the family of the Clark measures of the inner function $\thn^2$; we have $\nu_{1}^{n} = \mu_{2n}$. Combining formulas \eqref{eq83} and \eqref{eq42} in the proof of Theorem \ref{t2}, we obtain
$$\|f\|_{L^1(\T)} \le \sup_{|\beta| = 1}
\left(
\sqrt{\frac{\nu_{\beta}^{n}(2\Delta)}{\nu_{1}^{n}(\Delta)}} + 4\pi m(\Delta) \int_{\T\setminus 2\Delta}\frac{1}{|z - \xi_c|^2} \,d\nu_{\beta}^{n}(z)
\right).
$$  
Observe that $\nu_{1}^{n}(\Delta) \ge m(\Delta)$ and $\nu_{\beta}^{n}(2 \Delta) \le 3 m (\Delta)$. Let $\xi_1$, $\xi_2$ be the nearest points to $\xi_c$ in $a(\nu_{\beta}^{n})\setminus 2\Delta$. Then $|\xi_c - \xi_{1,2}| \ge \diam(2\Delta)/2 \ge m(2\Delta)$ and we have 
\begin{equation}\notag
\begin{aligned}
\int_{\T\setminus 2\Delta}\frac{d\nu_\beta^n(z)}{|z - \xi_c|^2}
&\le \int_{\T\setminus 2\Delta}\frac{dm(z)}{|z - \xi_c|^2} + \frac{1}{2n}\left(\frac{1}{|\xi_c - \xi_{1}|^2} + \frac{1}{|\xi_c - \xi_{2}|^2}\right)\\
&\le \frac{\pi}{4m(2\Delta)} + \frac{1}{2m(2\Delta)} < \frac{1}{m(\Delta)}.
\end{aligned}
\end{equation}
Hence, $\|f\|_{L^1(\T)} \le \sqrt{3} + 4 \pi < 15$. This gives us $d'_{n} < 15$ and $c_{1,n}^{-1} < 90$.

\medskip

\noindent Let us turn to the second inequality in \eqref{eq81}. As before,  from formula \eqref{eq84} we obtain 
\begin{equation}
\begin{aligned}
\|\Gamma_\phi\| &= \sup\left\{|\Psi(h)|, \; \|h\|_{Y_{\thn}} \le 1 \right\} 
\le D''_{n} \sup\left\{|\Psi(h)|, \; \|h\|_{L^1(\T)} \le 1\right\} \\
&\le D'_{n} D''_{n} \sup\left\{|\Psi(h)|, \; \|h\|_{H^{1}_{at}(\mu_{2n})} \le 1\right\} 
\le D'_{n} D''_{n} \|\phi\|_{\mu_{2n}^{*}}, 
\end{aligned}
\end{equation}
where 
\begin{equation}
\begin{aligned}
&D'_{n} = \sup\{\|h\|_{H^{1}_{at}(\mu_{2n})},\; h \in \Kthns \cap z H^1, \; \|h\|_{L^1(\T)} \le 1\}; \\
&D''_{n} = \sup\{\|h\|_{L^1(\T)},\; h \in \Kthns \cap z H^1, \; \|h\|_{Y_{\thn}} \le 1\}.
\end{aligned}
\end{equation}
By the Cauchy-Schwarz inequality, $D''_{n} \le 1$.  In the proof of Theorem \ref{t2prime} we have seen that $D'_{n} \le 24 c_{5n} c_{6n} M$, where $M$ is the norm of the non-tangential maximal operator $F \mapsto F^*$ on $H^1$ and $c_{5n}$, $c_{6n}$ are the constants $c_5$, $c_6$ from Lemma \ref{l8} for the inner function $\theta = \thn$. Since $\tilde A_{\mu_{2n}} = \tilde B_{\mu_{2n}} = 1$, we have $D'_{n} \le 24 \cdot 4 \cdot 60 \cdot M \cdot \eps_{n}^{-1}$, where $\eps_n$ stands for the parameter $\eps$ in Lemma~\ref{l7} for $\theta = \thn$. Next, since the sublevel set $\Omega_\delta$ of $\thn$ is connected for every $\delta > 0$, the proof of Lemma~\ref{l7} shows that one can take $\eps_{n} = \kappa_n/2$, where $\kappa_n \le \kappa^*_n =(2\tilde B_{\mu_{2n}})^{-1} = 1/2$ is chosen so that estimate~\eqref{eq46} holds for $\theta = \thn^{2}$, $\kappa = \kappa_n$. It remains to show that $\inf_{n}\kappa_n > 0$. For this aim it is sufficient to prove that the functions 
$f_{\xi_0, n} = f_{\xi_0}$ in formula~\eqref{eq59} for $\theta = \theta_n$ are bounded uniformly in  $n$. By formula \eqref{eq45}, $C_{\mu_{2n}} \le 1/2$. Next, for every pair of atoms $\xi, \xi_0 \in a(\mu_{2n})$ and for all $z \in D_{\xi_0}(\kappa_n^*)$ we have
$|\xi - z| \ge |\xi - \xi_0|/2$. Since $\inf_n A_{\mu_{2n}} = A_{\mu_{4}} = \frac{1}{4\sqrt{2}}$ and $B_{\mu_{2n}} \le \tilde B_{\mu_{2n}} = 1$ we see that estimate~\eqref{eq39} for $\sa = \mu_{2n}$ takes the following form:
\begin{equation}\notag
\begin{aligned}
\int_{\T \setminus \{\xi_0\}} \frac{\,d\mu_{2n}(\xi)}{|\xi - \xi_0|^2} 
&\le 
2\pi \int_{\T \setminus \Delta} \frac{\,d\mu_{2n}(\xi)}{|\xi - \xi_0|^2} + \frac{64}{\mu_{2n}\{\xi_0\}} 
\\
&\le \frac{\pi^2}{2 m(\Delta)} + \frac{64}{\mu_{2n}\{\xi_0\}}
\\
&\le \left(64 + \frac{\pi^2}{4}\right)\frac{1}{\mu_{2n}\{\xi_0\}}.
\end{aligned}
\end{equation}
It follows that $|f_{\xi_0}| < 70$ on $D_{\xi_0}$ and estimate \eqref{eq46} holds for $\theta = \thn$ with any constant $\kappa_n \le \kappa_n^*$ such that $70 \le (2 \kappa_n)^{-1}$. In particular, one can take $\kappa_n = 1/140$ for all $n \ge 2$. We now see that the constants $c_{2, n}$ are bounded: $c_{2, n} \le D'_{n} \le 24 \cdot 4 \cdot 60 \cdot 280 \cdot M < 10^7 M$. \qed

\medskip

\begin{Cor}\label{Cor2}
Let $A = (\alpha_{j-k})_{0 \le k, j \le n-1}$ be a Toeplitz matrix of size $n \times n$; consider its standard symbol $\psi = \alpha_{-(n-1)}\bar z^{n-1} + \ldots + \alpha_{n-1}z^{n-1}$. We have 
\begin{equation}\notag
c_1 \|\bar z^{n}\psi\|_{\mu_{2n}^{*}} \le \|A\| \le c_2 \|\bar z^{n}\psi\|_{\mu_{2n}^{*}},
\end{equation}
where the constants $c_1, c_2$ do not depend on $n$. 
\end{Cor}

\medskip

The author failed to find a simple argument allowing obtain Corollary~\ref{Cor1} from the $\bmo$-criterium for the boundedness of Hankel operators on $H^2$. The inverse implication is quite elementary.
\begin{Prop}\label{p3}
Let $\phi \in \ov{zH^2}$. The Hankel operator $H_\phi: H^2 \to \ov{z H^2}$ is bounded if and only if $\phi \in \bmo(\T)$. Moreover we have $c_1 \|\phi\|_{*} \le \|H_\phi\| \le c_2 \|\phi\|_{*}$ with constants $c_1, c_2$ from Corollary \ref{Cor1}.  
\end{Prop}
\beginpf Let $H_\phi: H^2 \to \ov{z H^2}$ be a bounded Hankel operator on $H^2$ with symbol $\phi \in \ov{zH^2}$. 
Then there are finite-rank Hankel operators $H_{\phi_n}$, $\phi_n \in \ov{\Kthn \cap z H^2}$, such that $H_\phi$ is the limit of $H_{\phi_n}$ in the weak$^*$ operator topology. Moreover one can choose $H_{\phi_n}$ so that $\sup_n\|H_{\phi_n}\| \le \|H_\phi\|$.  For every $n \ge 1$ and $k \ge n$ the operator norm of the Hankel operator $H_{\phi_n}$ is equal to the operator norm of the truncated Hankel operator on $K^{2}_{\theta_k}$ with  symbol $\phi_n$, where $\theta_k = z^k$. Since $\|\phi_n\|_{*} = \lim_{k \to \infty}\|\phi_n\|_{\mu_{2k}^{*}}$ we see from Corollary~\ref{Cor1} that 
\begin{equation}\label{eq90}
c_1 \|\phi_n\|_{*} \le \|H_{\phi_n}\| \le c_2 \|\phi_{n}\|_{*}, \quad n \ge 1.
\end{equation}
It follows that $c_1\sup \|\phi_n\|_{*} \le \|H_\phi\|$. Since $H_{\phi_n}$ tend to $H_{\phi}$ in the weak$^*$ operator topology we have $\lim_{n \to \infty}\int_\T p\phi_n\,dm = \int_\T p\phi\,dm$ for every trigonometric polynomial $p$. It is well-known that $H^{1}_{at}(\T)^* = \bmo(\T)$ (it worth be mentioned that this fact is much more easier than the Fefferman theorem on $\Re (zH^1)^* = \bmo(\T)$ which is generally used in the proof of the boundedness criterium for Hankel operators). Since trigonometric polynomials are dense in $\bmo(\T)$ in the weak$^*$ topology generated by $H^{1}_{at}(\T)$, we have $\phi \in \bmo(\T)$ and $c_1\|\phi\|_{*}\le \|H_\phi\|$. 
Now let $\phi \in \ov{z H^2} \cap \bmo(\T)$. Then there are functions $\phi_n \in 
\ov{\Kthn \cap z H^2}$ which tend to $\phi$ in the weak$^*$ topology of $\bmo(\T)$ and such that $\sup_{n}\|\phi_n\|_{*} \le \|\phi\|_{*}$. From \eqref{eq90} we see that $\|H_{\phi_n}\| \le c_2 \|\phi\|_{*}$ for the corresponding Hankel operators $H_{\phi_n}$. Since $L^2(\T) \subset H^{1}_{at}(\T)$ the functions $\phi_n$ converge to $\phi$ weakly in $L^2(\T)$. Hence for every pair of analytic polynomials $p_1, p_2$ we have 
$\lim_{n \to \infty} (H_{\phi_n}p_1, \ov{z p_2}) = (H_{\phi}p_1, \ov{z p_2})$. It follows that the operators $H_{\phi_n}$ converge  to the operator $H_\phi$ in the weak operator topology and we have $\|H_{\phi}\| \le c_2 \|\phi\|_{*}$. \qed

\bigskip

\bibliographystyle{plain} 
\bibliography{bibfile}

\enddocument